\documentclass[11pt,reqno]{amsart}

\usepackage{amssymb}
\usepackage{hyperref}

\theoremstyle{plain}
\newtheorem{thm}{Theorem}[section]
\newtheorem{lem}[thm]{Lemma}

\theoremstyle{definition}

\numberwithin{equation}{section}

\newcommand{\T}{\mathbb{T}}
\newcommand{\G}{\mathcal{G}}
\newcommand{\F}{\mathcal{F}}
\newcommand{\E}[1]{\mathbb{E}\left[#1\right]}
\newcommand{\pr}[1]{\mathbb{P}\left[#1\right]}
\newcommand{\Est}[1]{\mathbb{E}^*\left[#1\right]}
\newcommand{\eps}{\epsilon}
\newcommand{\lm}{\lambda}
\newcommand{\pgw}{\mathrm{PGW}_{\lambda}}
\newcommand{\poi}{\mathrm{Poisson}}

\title[Local algorithms and independent sets]{Local algorithms for independent sets are half-optimal}

\author{Mustazee Rahman \and B\'{a}lint Vir\'{a}g}

\address[Mustazee Rahman and B\'{a}lint Vir\'{a}g]{Department of Mathematics\\
University of Toronto\\
40 St. George Street\\
Toronto\\
ON M5S 2E4\\
Canada}

\email[Mustazee Rahman]{mustazee@math.toronto.edu}
\email[B\'{a}lint Vir\'{a}g]{balint@math.toronto.edu}

\thanks{M. Rahman's research was supported by a NSERC CGS grant.
B. Vir\'{a}g's research was partially supported by the Canada Research Chair program and the NSERC Discovery Accelerator Supplement.}

\begin{document}

\begin{abstract}
We show that the largest density of factor of i.i.d.~independent sets in the $d$-regular tree
is asymptotically at most $(\log d)/d$ as $d \to \infty$. This matches the lower bound given by previous constructions.

It follows that the largest independent sets given by local algorithms on random $d$-regular graphs have the same asymptotic
density. In contrast, the density of the largest independent sets in these graphs is asymptotically $2(\log d)/d$.

We prove analogous results for Poisson-Galton-Watson trees, which yield bounds for local algorithms
on sparse Erd\H{o}s-R\'{e}nyi graphs.
\end{abstract}

\maketitle

\section{Introduction} \label{sec:intro}

Local algorithms are randomized algorithms that run in parallel at each vertex of a graph by using only local information around each vertex. They produce important structures in large graphs, such as independent sets, matchings and colourings, with only constant running time (see \cite{Cso, CGHV, CL, EL, GG, HLS, LW, LN} and the references therein). In this paper we investigate local algorithms for high density independent sets in random $d$-regular graphs. We find an optimal bound for the density of such independent sets as the degree becomes large. It turns out that in this limit local algorithms can only yield independent sets with half the maximum possible density.

The motivation for our work comes from questions that arose in the theory of graph limits (see \cite{EL, HLS} and the references therein). In particular, Hatami, Lov\'{a}sz, and Szegedy conjecture (\cite{HLS} Conjecture 7.13) that most optimization problems over typical, sparse graphs can be solved by local algorithms.

We use the following notion of local algorithm introduced in
\cite{HLS}. The input to the algorithm is a graph $G$. The algorithm decorates $G$ by putting i.i.d.~labels on the vertices.  The output is $(f(i(v)); v \in G)$ where $f$ depends on the isomorphism class $i(v)$ of the labelled, rooted $r$-neighbourhood of $v$ for some fixed $r$. The process $(f(i(v)); v \in G)$ generated by the local algorithm will be called a \textbf{factor of i.i.d.}~process. See Section \ref{sec:regulartrees} for a more formal definition.

While the conjecture of Hatami, Lov\'{a}sz, and Szegedy was verified for maximal matchings \cite{BLS, CL,LN} and covariance structures
\cite{BSV}, Gamarnik and Sudan \cite{GS} showed that it fails for maximal independent sets. An independent set in a graph is a set of
vertices that have no edges between them.

It is known from \cite{BGT} that for each $d$ the size density of the largest independent sets in a random $d$-regular graph on $n$ vertices converges almost surely as $n \to \infty$. Furthermore, Bollob\'{a}s \cite{Bol} and McKay \cite{McKay} proved that with high probability the size density of the largest independent sets in random $d$-regular graphs is at most $2(\log d)/d$ for every $d \geq 3$. Frieze and {\L}uczak \cite{Fri, FL} provided lower bounds of matching asymptotic order for large $d$. Recently, precise formulae were given for large $d$ by Ding, Sly and Sun \cite{DSS}. On the other hand, several authors have produced local algorithms on $d$-regular graphs of large girth that yield independent sets of density $(\log d)/d$ for large $d$ (see \cite{GG, LW, Shearer}). These algorithms use greedy strategies to construct independent sets and can be easily adapted to random $d$-regular graphs.

Thus, for large $d$, the density of the largest independent sets in random $d$-regular graphs is of order $2(\log d)/d$ while
local algorithms have only produced independent sets with density of order $(\log d)/d$. The conjecture of Hatami, Lov\'{a}sz,
and Szegedy would imply that local algorithms can in fact produce independent sets in random $d$-regular graphs of density $2(\log d)/d$.

Gamarnik and Sudan \cite{GS} disprove this conjecture by showing that for large $d$ local algorithms can not find independent
sets in random $d$-regular graphs of density larger than $(1+ \frac{1}{\sqrt{2}})(\log d)/d$. Their crucial step is to prove that with high
probability any two high density independent sets in random $d$-regular graphs have a substantially large or substantially small intersection.
This observation was guided by predictions from statistical physics regarding the solution-space geometry of constraint satisfaction
problems \cite{MMbook}. In particular, the so called clustering phenomenon is expected to hold for independent sets in sparse random graphs.
Rigorous results have been established in this regard by Coja-Oghlan and Efthymiou \cite{COE} and in the aforementioned work
of Ding, Sly and Sun \cite{DSS}. It is shown that for large enough $d$, some of the properties that determine clustering
emerge for independent sets in random $d$-regular graphs at size density $(\log d)/d$.

In this paper we analyze the intersection densities of many independent sets in random regular graphs.
We show that with high probability (i.e., with probability tending to one as the size of the graphs tends to infinity)
the intersection densities must satisfy various inequalities. These structural results on the admissible intersection densities
imply quantitative bounds on the density of independent sets that can be generated from local algorithms. With the help of
these inequalities we prove that for any $\eps > 0$, local algorithms can not find independent sets in random $d$-regular
graphs of density larger than $(1+ \eps)(\log d)/d$ if $d$ is sufficiently large. In practice, iterative search algorithms that use local
moves at each step fail to find independent sets with density exceeding the critical threshold of $(\log d)/d$ in
random $d$-regular graphs. Our result provides some evidence as to why this is the case.

We also consider local algorithms for independent sets in Poisson-Galton-Watson trees. These yield local algorithms for independent sets
in sparse Erd\H{o}s-R\'{e}nyi graphs. We prove that the maximal density of local independent sets in a Poisson-Galton-Watson tree
of expected degree $\lm$ is of asymptotic order $(\log \lm)/\lm$ as $\lm \to \infty$. The aforementioned results of Bollob\'{a}s \cite{Bol},
Frieze and {\L}uczak \cite{Fri, FL} show that the largest independent sets in Erd\H{o}s-R\'{e}nyi graphs of average degree $\lm$ have
density of asymptotic order $2 (\log \lm)/\lm$ as $\lm \to \infty$.

The challenge in proving upper bounds to the density of local independent sets in Poisson-Galton-Watson trees is showing
that the randomness of the tree does not provide local algorithms with extra power. Also, in order to show the existence
of local independent sets having density close to $(\log \lm)/\lm$ we employ a coupling argument that produces
independent sets in Poisson-Galton-Watson trees from independent sets in regular trees.

\subsection{Organization of the paper}

In Section \ref{sec:regulartrees} we define the notion of a local algorithm for independent sets in the $d$-regular tree and relate it to local algorithms on finite $d$-regular graphs. Our main result about the density of local independent sets in regular trees is stated in Theorem \ref{thm:regularupbound}. In Section \ref{sec:upbound} we introduce the key inequality, stated in Theorem \ref{thm:binomineq}, that is satisfied by the intersection densities of any finite collection of local independent sets in the $d$-regular tree. Using this inequality we prove Theorem \ref{thm:regularupbound} in Section \ref{sec:proof}. In Section \ref{sec:binomineq} we prove Theorem \ref{thm:binomineq} by employing combinatorial arguments involving random regular graphs. In Section \ref{sec:pgwtrees} we state and prove our main result, Theorem \ref{thm:pgwthm}, on local independent sets in Poisson-Galton-Watson trees.

\section{Local algorithms for independent sets in regular graphs} \label{sec:regulartrees}

We define the notion of local algorithms for independent sets in regular trees.
Let $\T_d$ denote the rooted $d$-regular tree, and for $r \geq 0$ let $\T_{d,r}$ denote the rooted $r$-neighbourhood of $\T_d$.
A \textbf{labelling} of $\T_d$ is a vector $x \in [0,1]^{\T_d}$, and a \textbf{random labelling} is a labelling $X$ where
the co-ordinates $X(v), v \in \T_d$, are independent, uniformly distributed random variables on $[0,1]$. A \textbf{factor}
on $\T_d$ is a measurable function $f : [0,1]^{\T_d} \to \{0,1\}$ (w.r.t. the Borel $\sigma$-algebra) such that $f$ is invariant
under all root preserving automorphisms of $\T_d$. In other words, $f$ is spherically symmetric about the root.
We say that $f$ depends on the $r$-neighbourhood of the root if $f$ is defined on $[0,1]^{\T_{d,r}}$.

Any factor $f$ on $\T_d$ defines a set-valued stochastic process $I$ on $\T_d$ as follows. Any graph automorphism $\phi$ of $\T_d$
acts on labels $x \in [0,1]^{\T_d}$ by $\phi \cdot x (v) = x(\phi^{-1}(v))$. Since the automorphism group of $\T_d$ acts transitively on the vertices,
given any vertex $v$ let $\phi_v$ be an automorphism that maps $v$ to the root. For a random labelling $X$ of $\T_d$ define $I(v) = f(\phi \cdot X)$.
Due to $f$ being invariant under root preserving automorphisms $I$ is well defined. We call $I$ a factor of i.i.d.~process on $\T_d$.
Note that the distribution of $I$ is invariant under the action of the automorphism group of $\T_d$ (however, factor of i.i.d.~processes
are more restrictive than invariant process).

A factor of i.i.d.~independent set in $\T_d$ is a factor of i.i.d.~process $I$ such that $I$ is an independent in $\T_d$ with probability 1.
Since the distribution of $I(v)$ does not depend on the vertex $v$, we define the \emph{density} of $I$ as
$$density(I) = \pr{I(\text{root}) = 1} = \E{f(X)}\,.$$ It is easy to see that a factor that generates independent sets can be
approximated by similar factors that depend on finite size neighbourhoods of the root (see \cite[Section 12]{HLS}).
In this manner a factor of i.i.d.~independent set of density $\rho$ can be approximated by finite neighbourhood factor
of i.i.d.~independent sets whose densities converge to $\rho$. Hence, there is no harm in assuming that all our factors
for independent sets depend on finite size neighbourhoods of the root.

\paragraph{\textbf{Example: A construction of Lauer and Wormald}} In \cite{LW} the authors analyze the following algorithm that generates factor of i.i.d.~independent sets in $\T_d$. Fix $p \in (0,1)$ and an integer $k \geq 1$. Let $U_0 = V(\T_d)$ and for $1 \leq i \leq k$ do the following. Let $S_i \subset U_{i-1}$ be a random subset resulting from the output of a Bernoulli percolation on $U_{i-1}$ at density $p$. Set $U_i = U_{i-1} \setminus (S_i \cup N(S_i))$, where $N(S_i)$ is the one-neighbourhood of the set $S_i$ in $\T_d$. Consider the subset $I' = \cup_{i=1}^k S_i$. $I'$ may not be an independent set only because some $S_i$ may contain both vertices along an edge. If a vertex $v \in I'$ has one of its neighbours also included in $I'$ then exclude $v$ from $I'$. This results in an independent set $I \subset I'$.

The random set $I$ is a factor of i.i.d.~independent set since the decision rule to include a vertex is (deterministically) invariant of the vertex, and the rule depends on the outcome of $k$ independent Bernoulli percolations on $\T_d$. Furthermore, a little thought shows that the factor for $I$ depends only on the $(k+1)$-neighbourhood of a vertex.

Lauer and Wormald show that taking $k = \frac{c}{p}$ and then letting $p \to 0$, followed by $c \to \infty$,
results in independent sets whose densities converge to $\beta(d) := \frac{1 - (d-1)^{-2/(d-2)}}{2}$.
A simple analysis shows that $\frac{\log(d-1)}{d-2} - 2 (\frac{\log(d-1)}{d-2})^2 \leq \beta(d) \leq \frac{\log(d-1)}{d-2}$.

\paragraph{\textbf{From trees to finite graphs}}
Given a factor of i.i.d.~independent set $I$ in $\T_d$, we can construct a (random) independent set in any $d$-regular graph $G$ on $n$ vertices via the following procedure. Recall that $I$ uses a factor $f$ that computes $I(v)$ by looking only at the isomorphism class of the labelled $r$-neighbourhood of $v$ in $\T_d$. We begin with a random labelling $X$ of the vertices of $G$. Given any vertex $v \in G$ if its $(r+1)$-neighbourhood, $N_{r+1}(G,v)$, is a tree then set $I_G(v) = f(X(u); u \in N_r(G,v))$. This is allowed since $N_r(G,v) = \T_{d,r}$ by assumption. Otherwise, set $I_G(v) = 0$.

We verify that $I_G$ is an independent set. For any edge $(u,v)$ such that both $N_{r+1}(G,u)$ and $N_{r+1}(G,v)$ are trees, the pair of values $(I_G(u), I_G(v))$ is the same as the values $(I(a), I(b))$ for any edge $(a,b)$ of $\T_d$ with labels $(X_i; i \in N_r(G,u) \cup N_r(G,v))$ lifted to $N_r(\T_d,a) \cup N_r(\T_d,b)$ in the natural way. Thus, $(I_G(u), I_G(v)) \neq (1,1)$ as required. On the other hand, if one of $N_{r+1}(G,u)$ or $N_{r+1}(G,v)$ is not
a tree then at least one of $I_G(u)$ or $I_G(v)$ is 0. Consequently, $I_G$ is an independent set in $G$. Notice also that if $B(G)$ is the number of vertices of $G$ whose $(r+1)$-neighbourhood is not a tree then the expected size density of $I_G$ is $\E{|I_G|/n} = density(I)(1 - B(G)/n)$.

We are going to use this technique to project factor of i.i.d.~independent sets from $\T_d$ to finite, $d$-regular graphs. The resulting processes on the finite graphs will be referred to as independent sets from local algorithms. We are now prepared to state our main result for independent sets in $d$-regular graphs.
Define $\alpha_d$ as follows:
\begin{equation} \label{eqn:alphad}
\alpha_d \frac{\log d}{d} = \sup \{ \, density(I) : I \;\text{is a factor of i.i.d.~independent set in}\; \T_d \}\,.
\end{equation}

\begin{thm} \label{thm:regularupbound}
The following ineqaulity holds for $\alpha_d$:
\[ \limsup_{d \to \infty} \, \alpha_d \leq 1\,.\]
In other words, for any $\eps > 0$ there exists a $D$ such that if $d > D$ then there are no local algorithms that generate independent
sets in $\T_d$ having density larger than $(1+ \eps) \frac{\log d}{d}$.
\end{thm}

\subsection{Key inequality for intersection densities of local independent sets} \label{sec:upbound}

We prove Theorem \ref{thm:regularupbound} by way of contradiction. Assuming otherwise, we pass to a subsequence in $d$ and assume
that for some $\alpha > 1$ we have $\alpha_d > \alpha$ for every $d$ along the subsequence. It follows that
for each such $d$ there exists a factor of i.i.d.~independent set in $\T_d$, say $I_d$, such that the density of $I_d$
is $\alpha \frac{\log d}{d}$. Let $f_d : [0,1]^{\T_d} \to \{ 0, 1\}$ denote the factor associated to $I_d$. Recall we
may assume that $f_d$ depends on a finite size neighbourhood of the root. So we assume that $f_d$ depends on
the $r_d$-neighbourhood of $\T_d$.

Now we construct many copies of $I_d$ that are correlated with each other via a parameter that we will control.
Fix $p \in [0,1]$ and let $S_d = S_d(p)$ denote a random subset of the vertices of $\T_d$ generated by a Bernoulli percolation of density $p$.
Also, let $X_i$ for $i \geq 0$ denote independent random labellings of $\T_d$. We construct independent sets $I_{d,i}$ for $i \geq 0$
by letting $I_{d,i}$ be generated from the factor $f_d$ with labels $X_0(v)$ for $v \in V(\T_d) \setminus S_d$ and $X_i(v)$ for $v \in S_d$.

As $f_d$ is defined on $[0,1]^{\T_{d,r_d}}$ it follows that $I_{d,i}(\rm{root})$ depends only on the labels
$X_0(v)$ and $X_i(v)$ for $v \in V(\T_{d,r_d})$ and the subset $S_{d,r_d} = S_d \cap V(\T_{d,r_d})$.
Also, the joint distribution of the sets $I_{d,i}$ is exchangeable over $i$ and each $I_{d,i}$
follows the distribution of $I_d$. This implies that the intersection of any $k$ of these local independent sets have a common density,
which we denote $\alpha_{k,d,p} \,\frac{\log d}{d}$. Note that $\alpha_{1,d} = \alpha$.
To reduce notational clutter we will denote $r_d$ by $r$ and $S_{d, r_d}$ by $S$ until the end of Section \ref{sec:regulartrees}.

We will achieve a contradiction by first showing that these intersection densities are constrained to satisfy an inequality for each $k$.
Secondly, we will violate these inequalities by tuning the coupling parameter $p$ (under the assumption that $\alpha > 1$).
The next theorem introduces these key inequalities. Their proof, discussed in Section \ref{sec:binomineq}, is based on a structure
theorem about independent sets in random $d$-regular graphs.

\begin{thm} \label{thm:binomineq}
For each $k \geq 1$ the quantities $\alpha_{i,d,p}$ for $1 \leq i \leq k$ satisfy the following
\begin{equation} \label{binomeq}
\liminf_{d \to \infty} \; \inf_{p \in [0,1]}\; \sum_{i=1}^k (-1)^{i-1}\binom{k}{i} \alpha_{i,d,p}(2-\alpha_{i,d,p}) \geq 0\,.\end{equation}
\end{thm}

Theorem \ref{thm:binomineq} is proved by counting the expected number of $k$-tuples of
independent sets $(I_i, \ldots, I_k)$ in random $d$-regular graphs such that their intersection
densities are close to the quantities $\alpha_{i,d,p} \frac{\log d}{d}$ for $1 \leq i \leq k$.
We show that if (\ref{binomeq}) fails then the probability of observing such $k$-tuples
of independent sets in random $d$-regular graphs is vanishingly small as the size of the
graphs tend to infinity. On the other hand, Lemma \ref{lem:concentration} implies that the existence
of the local independent sets $(I_{d,1}, \ldots, I_{d,k})$ allows us to observe such $k$-tuples of independent sets
in random $d$-regular graphs with high probability and so (\ref{binomeq}) must hold.

\paragraph{\textbf{Relation to the approach of Gamarnik and Sudan}}

In their paper \cite{GS} Gamarnik and Sudan derive inequality (\ref{binomeq}) for $k= 2$. The $k=2$ case gives
$$\inf_{ p \in [0,1]} 2 \alpha(2-\alpha) -  \alpha_{2,d,p}(2-\alpha_{2,d,p}) \geq 0 \quad \text{for all large}\; d.$$
To minimize this in $p$ we certainly want to set $\alpha_{2,d,p} = 1$ for every $d$. It turns out that $\alpha_{2,d,p}$ is
continuous in $p$ (see Lemma \ref{lem:continuity}) with $\alpha_{2,d,0}  = \alpha$ and $\alpha_{2,d,1} = \alpha^2 (\frac{\log d}{d})$.
So if $\alpha > 1$ then for all large $d$ we can find a value of $p$ such that $\alpha_{2,d,p} = 1$.
This implies that the density $\alpha$ satisfies $\alpha(2 - \alpha) \geq 1/2$, or equivalently, that
$\alpha \leq 1 + \frac{1}{\sqrt{2}}$. This is the conclusion of Gamarnik and Sudan.

We may also analyze (\ref{binomeq}) for $k=3$ to conclude that $\alpha \leq 1 + \frac{1}{\sqrt{3}}$. Indeed, we have
that $3 \alpha(2-\alpha) - 2\alpha_{2,d,p}(2-\alpha_{2,d,p}) + \alpha_{3,d,p}(2-\alpha_{3,d,p}) \geq 0$ for large $d$.
If $\alpha >1$ then for all large $d$ we may choose a value of $p$ such that $\alpha_{2,d,p} = 1$. Also,
observe that $\alpha_{3,d,p}(2-\alpha_{3,d,p}) \leq 1$. Thus, we conclude from (\ref{binomeq}) that
$ 3 \alpha(2-\alpha) - 2 + 1 = 3 \alpha(2-\alpha) - 1 \geq 0$. This implies that $\alpha \leq 1 + \frac{1}{\sqrt{3}}$.

We do not know how to solve the minimization problem in $p$ exactly for $k \geq 4$.
In order to analyze (\ref{binomeq}) for large values of $k$ we are going to make
a choice of $p$ for each $d$ (and fixed $k$) that allows us to bound the
sum in (\ref{binomeq}) from above as $d \to \infty$. This upper bound is going to be a quantity that
we can analyze in the large $k$ limit. From there we will derive a contradiction to the assumption
that $\alpha > 1$.

\subsection{Proof of Theorem \ref{thm:regularupbound} from Theorem \ref{thm:binomineq}} \label{sec:proof}

Given the setup thus far we begin by interpreting the $\alpha_{k,d,p}$ in a probabilistic manner.
We show that the values $\frac{\alpha_{k,d,p}}{\alpha_{1,d}}$ can be realized as the moments of a random variable.
This random variable is defined on a new probability space, which is obtained from the
original probability space by essentially restricting
to the support of the factor $f_d$. Formally, the new sample space is the set $\{ f_d(X_0) \equiv 1\}$
considered as a subset of the joint sample space of $X_0, X_1, \ldots$, and $S$.
The new $\sigma$-algebra is the restriction of the $\sigma$-algebra generated by $S, X_0, X_1,\ldots$
to $\{f_d(X_0) \equiv 1\}$. The new expectation operator $\mathbb{E}^*$ is defined by
\[ \Est{U} = \frac{\E{f_d(X_0) \,U}}{\E{f_d(X_0)}} \]
for any random variable $U$ defined on $\{f_d(X_0) \equiv 1\}$.

If $\mathcal{F}$ is a $\sigma$-algebra such that $f_d(X_0)$ is $\mathcal{F}$-measurable,
then for any random variable $U$ defined on the original probability space we have
\[ \Est{U \mid \F} = \E{U \mid \F}\,. \]
This is to be interpreted by restricting $\F$ to $\{f_d(X_0) \equiv 1\}$ on the left hand side and the random variable
$\E{U \mid \F}$ to $\{f_d(X_0) \equiv 1\}$ on the right hand side. To prove this suppose that $Z$ is a $\F$-measurable random variable.
Then,

\begin{eqnarray*}
\Est{Z\,\E{U \mid \F}} &=& \frac{\E{f_d(X_0) Z\,\E{U \mid \F}}}{\E{f_d(X_0)}}\\
&=& \frac{\E{\E{f_d(X_0) Z U \mid \F}}}{\E{f_d(X_0)}} \quad (\text{since}\; f_d(X_0) \;\text{and}\; Z \;\text{are}\; \F \text{--measurable}) \\
&=& \frac{\E{f_d(X_0) Z U}}{\E{f_d(X_0)}}\\
&=& \Est{ZU}\,.
\end{eqnarray*}

Define a sequence of $[0,1]$-valued random variables $Q_{d,p} = Q_d(S, X_0)$, which we denote the \textbf{stability}, on the restricted probability space as follows.
Let $$f_{d,i} = f_d(X_0(v); v \notin S,\, X_i(v); v \in S) = I_{d,i}(\rm{root}).$$
Set $$Q_{d,p} = \Est{f_{d,1} \mid X_0, S} = \E{f_{d,1} \mid X_0, S}\,.$$
Roughly speaking, the stability is the conditional probability, given the root is included in the independent set,
that it remains to be included after re-randomizing the labels on $S$.

The key observation is that the moments of the stability satisfy $\Est{Q_{d,p}^{k-1}} = \frac{\alpha_{k,d,p}}{\alpha_{1,d}}$ for $k \geq 1$.
Indeed, as $Q_{d,p}$ has the same distribution as $\E{f_{d,i} \mid X_0, S}$ for every $i$ we have that
$$ \Est{Q_{d,p}^{k-1}} = \frac{\E{f_{d,0} \, (\E{f_{d,1} \mid X_0, S})^{k-1}}}{\E{f_{d,0}}} =
\frac{\E{f_{d,0} \, \left(\prod_{i=1}^{k-1}\E{f_{d,i} \mid X_0, S}\right)}}{\E{f_{d,0}}}\,.$$

The random variables $f_{d,i}$ are independent of each other conditioned on $(X_0,S)$. Hence,
$$\prod_{i=1}^{k-1}\E{f_{d,i} \mid X_0, S} = \E{\prod_{i=1}^{k-1} f_{d,i} \mid X_0, S}\,.$$
Furthermore, $f_{d,0}$ is measurable w.r.t.~$(X_0,S)$ and so we conclude that
$$\E{f_{d,0} \, \left(\prod_{i=1}^{k-1}\E{f_{d,i} \mid X_0, S}\right)} = \E{f_{d,0} \prod_{i=1}^{k-1} f_{d,i}} = \text{density}(\cap_{i=1}^k I_{d,i})\,.$$
Consequently, $\Est{Q_{d,p}^{k-1}} = \frac{\text{density}(\cap_{i=1}^k I_{d,i})}{\text{density}(I_{d,1})} = \frac{\alpha_{k,d,p}}{\alpha_{1,d}}$.

Henceforth, all expectations involving $Q_{d,p}$ will simply be denoted by $\mathbb{E}$ instead of $\mathbb{E}^*$.
We will need the following lemma regarding the continuity of the stability in terms of the coupling parameter $p$.

\begin{lem} \label{lem:continuity}
Let $g : [0,1] \to \mathbb{R}$ be a continuous function.
The moment $\E{g(Q_{d,p})}$ is a continuous function of $p$.
When $p=0$, $\E{g(Q_{d,0})} = g(1)$, and when $p=1$, $\E{g(Q_{d,1})} = g(\alpha \frac{\log d}{d})$.
\end{lem}

\begin{proof}
The parameter $p$ enters into $\E{g(Q_{d,p})}$ only through the random finite subset $S \subset T_{d,r}$. For each $W \subset \T_{d,r}$ the
probability $\pr{S = W} = p^{|W|}(1-p)^{|T_{d,r} \setminus W|}$. This probability is a polynomial in $p$. By conditioning on the output of $S$ we note that $\E{g(Q_{d,p})}$ can be expressed as a convex combination of terms that are free of $p$, namely $\E{g(Q_{d,p})| S = W}$,
with corresponding coefficient $\pr{S = W}$. Thus, $\E{g(Q_{d,p})}$ is also a polynomial in $p$.

When $p=0$ the set $S$ is empty and $f_{1,d} = f_{d,0}$. Therefore, conditioning on $X_0$ and restricting to $\{f_{d,0} \equiv 1\}$ forces $Q_{d,p} \equiv 1$. When $p=1$ the set $S$ equals $\T_{d,r}$, and hence $f_{1,d}$ becomes independent of the random labelling $X_0$, and hence of $f_{d,0}$ as well. Consequently, the conditioning has no effect and $Q_{d,p} = \E{f_d} = \alpha \frac{\log d}{d}$.
This implies that $\E{g(Q_{d,0})} = g(1)$ and $\E{g(Q_{d,1})} = g(\alpha \frac{\log d}{d})$.
\end{proof}

We now translate the inequality from (\ref{binomeq}) in terms of the stability.
Our goal is to rewrite (\ref{binomeq}) as an expectation of a function of the stability, which we can then
analyze for large values of $d$ and $k$.

Note that $\alpha_{i,d,p} = \alpha_{1,d} \E{Q_{d,p}^{i-1}} = \alpha \E{Q_{d,p}^{i-1}}$.
To deal with the terms $\alpha_{i,d,p}^2$ we introduce an independent copy of $Q_{d,p}$, which
we denote $R_{d,p}$. Thus, $\alpha_{i,d,p}^2 = \alpha^2 \E{(Q_{d,p}R_{d,p})^{i-1}}$.
This implies that $$\alpha_{i,d,p}(2-\alpha_{i,d,p}) = 2 \alpha \E{Q_{d,p}^{i-1}} - \alpha^2 \E{(Q_{d,p}R_{d,p})^{i-1}}.$$

Observe the following identity that results from the binomial theorem:
\begin{equation} \label{eqn:binomidn}
\sum_{i=1}^{k} (-1)^{i-1} \binom{k}{i} x^{i-1} = \frac{1 - (1-x)^k}{x} \quad \text{for}\; 0 \leq x \leq 1.
\end{equation}

Let $s_k(x) =  \frac{1 - (1-x)^k}{x}$ for $x \in [0,1]$ and $k \geq 1$. Note that $s_k(0) = \lim_{x \to 0} s_k(x) = k$.
We may now translate the inequality from (\ref{binomeq}) into
\begin{equation} \label{Qinq}
\liminf_{d \to \infty} \; \inf_{p \in [0,1]}\; 2\alpha \E{s_k(Q_{d,p})} - \alpha^2 \E{s_k(Q_{d,p}R_{d,p})} \geq 0\,.
\end{equation}

We make a particular choice of $p$ for every $d$ in order to analyze (\ref{Qinq}) in the large $d$ limit.
Fix a parameter $u > 0$ that we will tune later. In the statement of Lemma \ref{lem:continuity} take $g(x) = x^u$ for $0 \leq x \leq 1$.
From the assumption that $\alpha > 1$, we employ Lemma \ref{lem:continuity} and deduce that for all sufficiently large $d$
we can select a $p = p(d,u)$ such that $$\E{Q_{d,p(d,u)}^u} = 1/\alpha.$$ 

We denote $Q_{d,p(d,u)}$ by $Q_d$. At this point our reasoning behind this choice is mysterious.
The idea, of course, is that by choosing $p$ this way we try to minimize the left hand side of
(\ref{Qinq}) in a manner that we can analyze as $k \to \infty$. The argument that follows will show that
our choice is judicious.

Recall that probability distributions on $[0,1]$ are compact with respect to convergence in distribution.
Therefore, from the sequence $(Q_d, R_d)$ we can choose a subsequence $(Q_{d_i}, R_{d_i})$ that converges in
distribution to limiting random variables $(Q, R)$. The random variables $Q$ and $R$ are independent and identically distributed
with values in $[0,1]$.

Observe that $s_k(x) = 1 + (1-x) + \cdots + (1-x)^{k-1}$. Thus, $s_k(x)$ is a continuous, decreasing function on $[0,1]$
with maximum value $s_k(0) = k$ and minimum value $s_k(1) = 1$. Therefore, distributional convergence of $(Q_{d_i},R_{d_i})$
to $(Q,R)$ implies that $\E{s_k(Q_{d_i})} \to \E{s_k(Q)}$ and $\E{s_k(Q_{d_i}R_{d_i})} \to \E{s_k(QR)}$.

By passing to the subsequence $d_i$ and taking limits in $i$ the inequality (\ref{Qinq}) becomes
\begin{equation} \label{Qinq2} 2 \E{s_k(Q)} \geq \alpha \E{s_k(QR)} .\end{equation}

This holds for every $k \geq 1$. Taking the limit as $k \to \infty$ of (\ref{Qinq2}) results in the inequality $2 \E{1/Q} \geq \alpha \E{1/Q}^2$. If $\E{1/Q}$ is finite then we have $\alpha \leq 2 \E{1/Q}^{-1}$. We are thus left with the seemingly contradictory task of showing that $\E{1/Q}$ is finite but large. Unfortunately, distributional convergence of $Q_{d_i}$ to $Q$ is not sufficient to get a lower bound on $\E{1/Q}$. Furthermore, it is not a priori clear that this expectation is finite, or even that $\pr{Q=0} = 0$. To work around these difficulties we have to control the distribution of the $Q_d$ well enough to be able to conclude that $\pr{Q=0}$ is small while $\E{1/Q}$ is large. We derive a contradiction to the hypothesis $\alpha > 1$ from analyzing (\ref{Qinq2}) based upon 3 cases: $\pr{Q = 0} > 0$, or $\pr{Q=0} = 0$ but $\E{1/Q} = \infty$, or $\E{1/Q} < \infty$.

In order to bound $\E{1/Q}$ and $\pr{Q=0}$ we recall that we had set $p$ such that $\E{Q_d^u} = 1/ \alpha$ for all large $d$.
Since $x \to x^u$ is continuous and bounded on $[0,1]$ we conclude that $\E{Q^u} = \lim_{i} \E{Q_{d_i}^u} = 1/ \alpha$.
Now, an upper bound on $\E{Q^u}$ implies a lower bound on $\E{1/Q}$ due to $\E{1/Q} \geq \E{Q^u}^{-1/u}$,
which follows from Jensen's inequality. Also, a lower bound on $\E{Q^u}$ gives a upper bound on $\pr{Q=0}$
because $\mathbf{1}_{Q = 0} \leq 1 - Q^u$. We now analyze (\ref{Qinq2}) over all $k$ based upon the 3 cases mentioned
in the previous paragraph.

\paragraph{\textbf{Case 1:}} $\pr{Q = 0} = q > 0$. In this case most of the contribution to $\E{s_k(Q)}$ results from $\{Q=0\}$. More precisely,
$\frac{s_k(x)}{k} = \mathbf{1}_{x=0} + \frac{s_k(x)}{k}\mathbf{1}_{x>0}$, and $\frac{s_k(x)}{k}\mathbf{1}_{x>0} \to 0$ as $k \to \infty$. Also,
$\frac{s_k(x)}{k} \in [0,1]$ for all $k$ and $x \in [0,1]$. Therefore, from the bounded converge theorem we deduce that
$\E{s_k(Q)/k} \to \pr{Q=0}$ as $k \to \infty$, and similarly, $\E{s_k(QR)/k} \to \pr{QR=0}$. The latter probability is $2q - q^2$ due to $Q$ and $R$ being
independent and identically distributed. Upon dividing the inequality in (\ref{Qinq2}) through by $k$ and taking a limit we conclude that
$$2q - \alpha(2q-q^2) \geq 0, \; \text{or equivalently that} \; \alpha \leq \frac{2}{2-q}\,.$$

For $x \in [0,1]$ we have that $\mathbf{1}_{x=0} \leq 1 - x^u$. It follows from this that $q \leq 1 - \E{Q^u} = 1 - 1/\alpha$. Thus,
$$ \alpha \leq \frac{2}{2-q} \leq \frac{2}{1 + \alpha^{-1}}\,.$$
Simplifying the latter inequality gives $\alpha \leq 1$; a contradiction.

\paragraph{\textbf{Case 2:}} $\pr{Q=0} = 0$ but $\E{\frac{1}{Q}} = \infty$. In this case most of the contribution to $\E{s_k(Q)}$ occurs when $Q$ is small.
Note that $s_k(x) \nearrow 1/x$ as $k \to \infty$. Hence, the monotone convergence theorem implies that $\E{s_k(Q)} \to \infty$ as $k \to \infty$.

Fix $0 < \epsilon < 1$, and write $s_k(x) = s_{k, \leq \eps}(x) + s_{k, > \eps}(x)$ where $s_{k, \leq \eps}(x) = s_k(x)\,\mathbf{1}_{x \leq \eps}$. Note that $s_{k, > \eps}(x) \leq \eps^{-1}$ for all $k$. We have that
\begin{equation} \label{eqn:supbound} \E{s_k(Q)} = \E{s_{k, \leq \eps}(Q)} +  \E{s_{k, > \eps}(Q)} \leq \E{s_{k, \leq \eps}(Q)} + \eps^{-1}.\end{equation}
Thus, $\E{s_{k, \leq \eps}(Q)} \to \infty$ with $k$ because $\E{s_k(Q)} \to \infty$.

We also observe from the positivity of $s_k$ that $$\E{s_k(QR)} \geq \E{s_k(QR); Q \leq \eps, R > \eps} + \E{s_k(QR); Q > \eps, R \leq \eps}.$$ The latter two terms are equal by symmetry, so $\E{s_k(QR)} \geq 2 \E{s_k(QR); Q \leq \eps, R > \eps}$. The fact that $s_k(x)$ is decreasing in $x$ and $R \leq 1$ implies that $s_k(QR) \geq s_k(Q)$. Together with the independence of $Q$ and $R$ we deduce that
$$\E{s_k(QR); Q \leq \eps, R > \eps} \geq \E{s_k(Q); Q \leq \eps, R > \eps} = \E{s_{k, \leq \eps}(Q)} \pr{R > \eps}.$$
Consequently,
\begin{equation} \label{eqn:slwbound} \E{s_k(QR)}  \geq 2 \E{s_{k, \leq \eps}(Q)} \pr{Q > \eps}.\end{equation}

The inequality in (\ref{Qinq2}) is $\frac{\alpha}{2} \leq \frac{\E{s_k(Q)}}{\E{s_k(QR)}}$. The bounds from (\ref{eqn:supbound}) and (\ref{eqn:slwbound})
imply that
$$\frac{\alpha}{2} \leq \frac{\E{s_{k, \leq \eps}(Q)} + \eps^{-1}}{2 \E{s_{k, \leq \eps}(Q)} \pr{Q > \eps}}\,.$$
Since $\E{s_{k, \leq \eps}(Q)} \to \infty$ with $k$ we can take a limit in $k$ to conclude that
$$ \alpha \leq \frac{1}{\pr{Q > \eps}}.$$
As $\eps \to 0$ the probability $\pr{Q > \eps} \to \pr{Q > 0} = 1$, by assumption. Thus, $\alpha \leq 1$; a contradiction.

\paragraph{\textbf{Case 3:}} $\E{1/Q}$ is finite. In the final case we again use the fact that $s_k(x)$ increases to $1/x$ for $0 \leq x \leq 1$. Hence,
$s_k(Q) \nearrow 1/Q$ almost surely and $s_k(QR) \nearrow 1/(QR)$ almost surely. Taking a limit of the inequality in (\ref{Qinq2})
and using the monotone convergence theorem it follows that
$$ 2 \E{\frac{1}{Q}} - \alpha \E{\frac{1}{QR}} \geq 0\,.$$
Since $\E{\frac{1}{QR}} = \E{\frac{1}{Q}}^2$ the above inequality reduces to $$\alpha \leq 2 \E{\frac{1}{Q}}^{-1} \leq 2 \E{Q^u}^{1/u}\,.$$ In the last step we have used the power-mean/Jensen's inequality. Since $\E{Q^u} = 1/\alpha$, we see that
$$\alpha \leq 2^{\frac{u}{u+1}}\,.$$
Due to the contradiction resulting from the previous two cases we deduce that for all $u > 0$ we have $\alpha \leq 2^{\frac{u}{u+1}}$. By letting
$u \to 0$ we conclude that $\alpha \leq 1$; the final contradiction.

\section{Inequalities for intersection densities: proof of Theorem \ref{thm:binomineq}} \label{sec:binomineq}

We will prove Theoem \ref{thm:binomineq} by reducing it to a problem about densities of independent sets in large, finite, $d$-regular graphs.
First, we begin with some terminology. Let $\G_{n,d}$ denote a random $d$-regular graph on $n$ vertices sampled according to the configuration model (see \cite{Bolbook} chapter 2.4): each of the $n$ distinct vertices emit $d$ distinct half-edges, and we pair up these $nd$ half-edges uniformly at random. These $nd/2$ pairs of half-edges can be glued into full edges to yield a labelled, random, $d$-regular graph. Note that the resulting graph can have loops and multiple edges. There are $(nd-1)!! = (nd-1)(nd-3)\cdots 3\cdot 1$ possible pairings, or outcomes, of the model. Let $G_{n,d}$ denote the set of all these outcomes.
So $\G_{n,d}$ is picked uniformly at random from $G_{n,d}$.

\subsection{Projecting independent sets from $\T_d$ to $\G_{n,d}$} \label{sec:projection}

Recall that given the factor of i.i.d.~independent set $I_d$ on $\T_d$ we can project it to a (random) independent set $I_G$ on any given $G \in G_{n,d}$.
If $B(G)$ is the number of vertices of $G$ whose $(r+1)$-neighbourhood is not a tree then $\E{|I_G|/n} = density(I_d)(1 - \frac{B(G)}{n})$.
We can model the independent sets $I_{d,i}$ from Section \ref{sec:upbound} in the random graph $\G = \G_{n,d}$. To do so we first choose a random subset $S_\G \subset V(\G) = [n]$ via a Bernoulli percolation with density $p$. Then we fix independent random labellings $X_i$ of $\G$ for $i \geq 0$. We define $I_{\G,i}$ to be the projection of $I_d$ with input $X_0(v)$ for $v \notin S_\G$ and $X_i(v)$ for $v \in S_\G$. As the $I_{\G,i}$ are exchangeable, for any finite subset $T \subset \{1, 2, \ldots\}$ we have $\E{|\cap_{i \in T} I_{\G,i}|/n} = \alpha_{|T|,d,p} \frac{\log d}{d}(1 - \frac{\E{B(\G)}}{n})$. It is well-known that $\E{B(\G)}$ is bounded in $n$ for every $d$ \cite[chapter 9.2]{JLbook} . So $\E{|\cap_{i \in T} I_{\G,i}|/n} \to \alpha_{|T|,d,p} \frac{\log d}{d}$ as $n \to \infty$ for every $T \subset [k]$.

\subsection{The expected number of independent sets satisfying a given density profile} \label{sec:expectation}

From the construction above we see that the factor of i.i.d.~independent set $I_d$ on $\T_d$ can be used to produce $k$-tuples of independent sets
$(I_{\G,1}, \ldots, I_{\G,k})$ in $\G = \G_{n,d}$ such that the intersection densities of these $k$ independent sets are close to those of $I_{d,1},\ldots, I_{d,k}$, defined in Section \ref{sec:upbound}. We will compute the expected number of $k$-tuples of independent sets in $\G$ with some given intersection densities. This will allow us to bound, from above, the probability of observing  a $k$-tuple of independent sets in $\G$ whose intersection densities are close to that of $I_{d,1},\ldots, I_{d,k}$.  We will show that the expected number of $k$-tuples of independent sets in $\G$ with some prescribed intersection densities is an exponential term of the form $e^{nR}$ (Theorem \ref{thm:expectation}). The dominating contribution to the rate $R$ is from the binomial sum of the left hand side of inequality (\ref{binomeq}) for the prescribed intersection densities (Lemma \ref{lem:maxentropy} and Lemma \ref{lem:asymptotics}). This will allow us to conclude that the only $k$-tuples of independent sets in $\G$ that exists (with non-vanishing probability as $n \to \infty$) are those for which the corresponding rate $R$ is non-negative. Then in the final step we will show via concentration inequalities that there exists $k$-tuples of independent sets in $\G$ whose intersection densities are close to those given by $I_{d,1},\ldots, I_{d,k}$ (Lemma \ref{lem:concentration}). This is the strategy behind the proof of Theorem \ref{thm:binomineq}. Before proceeding we introduce some terminology.

For a $k$-tuple of independent sets $(I_{G,1}, \ldots, I_{G,k})$ in $G \in G_{n,d}$, the \textbf{density profile} associated to this $k$-tuple is the
vector $\rho = (\rho(T); T \subset [k])$ defined by $\rho(T) = |\cap_{i \in T} I_{G,i}|/n$ (set $\rho(\emptyset) = 1$). Associated to this $k$-tuple is also an ordered partition $\Pi$ of $V(G)$ into $2^k$ cells defined as follows:
$$\Pi = \{ \Pi(T): T \subset [k]\} \quad \text{with} \quad \Pi(T) = \left (\bigcap_{i \in T}  I_{G,i} \right ) \cap
\left ( \bigcap_{i \notin T} (V(G) \setminus I_{G,i}) \right ).$$
In other words, $\Pi(T)$ consists of vertices that belong to all the sets $I_{G,i}$ for $i \in T$ and none of the other sets. The partition $\Pi$ defines a probability measure $\pi = (\pi(T); T \subset [k])$ on $2^{[k]}$ by $\pi(T) = |\Pi(T)|/n$. This correspondence between $k$-tuples $(I_{G,1}, \ldots, I_{G,k})$ and ordered partitions $\Pi$ is bijective, and by the inclusion-exclusion principle we have that
\begin{eqnarray}
\label{eqn:PIEpi} \pi(T) &=& \sum_{T': T \subset T'} (-1)^{|T' \setminus T|} \rho(T') \, ,\\
\label{eqn:PIErho} \rho(T) &=& \sum_{T': T \subset T'} \pi(T')\,.
\end{eqnarray}

Finally, corresponding to $G$ and $\Pi$ is a $2^k \times 2^k$ matrix $M$ that we denote the \textbf{edge profile} of $\Pi$. For $T, T' \subset [k]$,
define $$M(T,T') = \frac{|\{(u,v) \in E(G): u \in \Pi(T), v \in \Pi(T')\}|}{nd}\,.$$

The tuple $(u,v)$ refers to a directed edge; so $(u,v) \neq (v,u)$ unless $u=v$. The number of directed edges of $G$ is $2|E(G)| = nd$.
Notice that $M(T,T')$ is the probability that a uniformly chosen directed edge of $G$ starts in $\Pi(T)$ and ends in $\Pi(T')$. Clearly, $M$ is a symmetric matrix with non-negative entries that sum to 1. Also, the marginal of $M$ along either the rows or columns is $\pi$. A crucial observation is that if $T \cap T' \neq \emptyset$ then $M(T,T') = 0$. Indeed, in this case both $\Pi(T)$ and $\Pi(T')$ lie in the common independent set $I_{G,i}$ for any $i \in T \cap T'$, and thus, there cannot be any edges joining $\Pi(T)$ to $\Pi(T')$.

Conversely, suppose we begin with an ordered partition $\Pi$ as above that induces an edge profile $M$ on $G$. If the edge profile satisfies the constraints $M(T, T') = 0$ whenever $T \cap T' \neq \emptyset$ then the $k$-tuple of subsets $(I_{G,1}, \ldots, I_{G,k})$ of $V(G)$ corresponding to $\Pi$ will be independents sets in $G$. Indeed, for any $i$, the number of edges of $G$ that have both endpoints in $I_{G,i}$ is
$(nd)/2 \sum_{(T,T'): i \in T \cap T'} M(T,T') = 0$. In this case the density profile $\rho$ of $(I_{G,1}, \ldots, I_{G,k})$ is given by (\ref{eqn:PIErho}) with $\pi$ being the marginal of $M$ along its rows.

With this terminology and bijection in mind let $Z(\rho) = Z(\G, \rho)$ denote the number of $k$-tuples of independent sets in $\G$ with density profile $\rho$.
Let $Z(\rho, M)$ denote the number of ordered partitions of $\G$ into $2^k$ cells such that the partitions induce the edge profile $M$, and $M$ is compatible with $\rho$ in the following sense.
The marginal, $\pi$, of $M$ along its rows is given by $\rho$ via (\ref{eqn:PIEpi}), and $M(T,T') = 0$ whenever $T \cap T' \neq \emptyset$.
It is clear from the discussion above that $$Z(\rho) = \sum_{M} Z(\rho, M)$$ where the sum is over all $M$ that is compatible with $\rho$.

\begin{thm} \label{thm:expectation}
Given the setup as above, define the entropies
\[ H(M) = \sum_{(T,T')} -M(T,T') \log(M(T,T')) \;\; \text{and} \;\;
H(\pi) = \sum_{T} -\pi(T) \log(\pi(T)) \quad (0\log 0 = 0).\]
The expectation of $Z(\rho,M)$ satisfies
\begin{equation} \label{eqn:expectation}
\E{Z(\rho,M)} \leq \rm{poly}(n,d,M_{\min}) \, \exp{\left\{ n \left[\frac{d}{2}H(M) - (d-1)H(\pi)\right] \right\}}\,.
\end{equation}
The term $\rm{poly}(n,d,M_{\min})$ is a polynomial in $n,d$, and $\frac{1}{M_{\min}}$ where $M_{\min} = \min \{M(T,T') : M(T,T') > 0\}$. The degree of this polynomial is bounded by a function of $k$ (at most $4^k$).
\end{thm}

\begin{proof} To compute the expectation we sum the probabilities of outcomes where each outcome uniquely specifies a pairing of half-edges in the configuration model that gives rise to a partition $\Pi$ with edge profile $M$. To specify such an outcome, do the following.
\begin{enumerate}
\item Partition the vertex set $[n]$ into $2^k$ distinguishable cells $\Pi(T),\; T \subset [k]$ with $| \Pi(T)| = n \pi(T)$.
\item Given the partition $\Pi$ from (1), and each subset $T \subset [k]$, partition the $nd\pi(T)$ half-edges attached to the vertices of $\Pi(T)$ into $2^k$ distinguishable cells
$\Pi(T,T'), T' \subset [k],$ such that $|\Pi(T,T')| = nd M(T,T')$.
\item For each pair $\{ T, T'\}$ with $T \neq T'$ pair up the half-edges from $\Pi(T,T')$ with those from $\Pi(T',T)$ in a specific way. Then for each $T$ pair the half-edges from $\Pi(T,T)$ with themselves in a specific way.
\end{enumerate}

Each outcome has probability $1/(nd-1)!!$ from definition of the configuration model. We compute the number of outcomes in the following.
But first, we should mention some conventions that we use in the following calculations. For an even integer $m \geq 2$ we denote
$(m-1)!! = (m-1)(m-3)\cdots 1$, and if $m=0$ then $(m-1)!! = 1$. Also, note that in any valid edge profile $M$ the quantities $ndM(T,T')$
have to be non-negative integers. Furthermore, $ndM(T,T)$ has to be even for every $T$ because for any $G \in G_{n,d}$
the number of half edges from $\Pi(T)$ to itself is twice the number of edges present in the subgraph of $G$ induced by $\Pi(T)$.
We may assume that $M$ has all these properties. We now compute the number of outcomes.

\begin{itemize}
\item The number of partitions of $[n]$ that satisfies the properties in (1) above is the multinomial coefficient
$$ \binom{n}{n\pi(T); T \subset [k]}\,.$$

\item Given a partition $\Pi$ satisfying (1) from above, the number of partitions of the half-edges that satisfy the properties in (2) is
$$ \prod_{T \subset [k]} \binom{nd\pi(T)}{ndM(T,T'); T' \subset [k]}\,.$$

\item Given the two partitions arising from (1) and (2), the number of pairings that satisfy (3) is
$$ \left[\prod_{(T,T'): T \neq T'} (ndM(T,T'))!\right]^{1/2} \prod_{T \subset [k]} (ndM(T,T)-1)!!\,.$$
\end{itemize}

The total number of outcomes is the product of the three terms above. From the linearity of expectation we conclude that
$\E{Z(\rho,M)}$ equals
\begin{eqnarray*}
\binom{n}{n\pi(T); T \subset [k]} &\times& \prod_{T \subset [k]} \binom{nd\pi(T)}{ndM(T,T'); T' \subset [k]} \times \\
&&\left[\prod_{(T,T'): T \neq T'} (ndM(T,T'))!\right]^{1/2} \times \prod_{T \subset [k]} (ndM(T,T)-1)!! \times \frac{1}{(nd-1)!!}\,.
\end{eqnarray*}

Now we do the asymptotics in $n$ by using Stirling's approximation of $m! \sim \sqrt{2\pi m} (m/e)^m$.
More precisely, $\sqrt{2\pi m} (m/e)^m \leq m! \leq (1 + \frac{1}{12m})\sqrt{2\pi m} (m/e)^m$.
Also, for an even integer $m$, $(m-1)!! = \frac{m!}{2^{m/2}(m/2)!}$.
In the following we need to consider only those values of $\pi(T)$ and $M(T,T')$ that are strictly positive.
We begin by simplifying the term

\begin{align*}
&\prod_{T \subset [k]} \binom{nd\pi(T)}{ndM(T,T'); T' \subset [k]} \left[\prod_{(T,T'): T \neq T'} (ndM(T,T'))!\right]^{1/2} \prod_{T \subset [k]} (ndM(T,T)-1)!!\\
&= \prod_{T} (nd\pi(T))! \left[\prod_{(T,T'): T \neq T'} (ndM(T,T'))!\right]^{-1/2} 2^{-nd/2 \sum_{T} M(T,T)} \left[\prod_{T} (\frac{nd}{2}M(T,T))!\right]^{-1}\,.\end{align*}

After incorporating the remaining two terms we see that the expectation is
\begin{align*} \binom{n}{n\pi(T); T \subset[k]} \times \binom{nd}{nd\pi(T); T \subset [k]}^{-1} \times (nd/2)! \times 2^{\left(\frac{nd}{2}(1-\sum_T M(T,T))\right)} &\times \\
\left[\prod_{(T,T'): T \neq T'} (ndM(T,T'))!\right]^{-1/2} \times \left[\prod_{T} (\frac{nd}{2}M(T,T))!\right]^{-1} \,.\end{align*}

Using Stirling's approximation we can verify that (with universal constants)
\begin{align*} \binom{n}{n\pi(T); T \subset[k]} \binom{nd}{nd\pi(T); T \subset [k]}^{-1} &= O(d^{(2^k-1)/2}) \prod_T \pi(T)^{n\pi(T)(d-1)}\\
&= O(d^{(2^k-1)/2}) \exp{\left\{-n(d-1)H(\pi)\right\}}\,.\end{align*}
Similarly,
$$\prod_T (\frac{nd}{2}M(T,T))! = O\left((nd\pi)^{2^{k-1}} \prod_T M(T,T)^{1/2}\right) \prod_T \left(\frac{nd M(T,T)}{2e}\right)^{\frac{nd}{2}M(T,T)}\,;$$

\begin{align*}
\prod_{\substack{(T,T')\\ T \neq T'}} (ndM(T,T'))! = O\left((2\pi nd)^{\binom{2^k}{2}} \prod_{\substack{(T,T')\\ T \neq T'}} M(T,T')^{1/2}\right) &\times (\frac{nd}{e})^{\{nd \sum_{\substack{(T,T')\\ T \neq T'}} M(T,T')\}}\times \\
&\prod_{\substack{(T,T')\\ T \neq T'}} M(T,T')^{ndM(T,T')}\,.
\end{align*}

Now, $(nd/2)! = O(\pi nd) (\frac{nd}{2e})^{nd/2}$. From this and the previous two equations we can check that all terms involving powers of $\frac{nd}{e}$ and powers of 2 algebraically cancel out from the expression for $\E{Z(\rho,M)}$. Therefore, after algebraic simplifications we conclude that
\begin{align*}
\E{Z(\rho,M)} &= O\left((\pi n)^{(1-2^k)/2} (2\pi nd)^{-\frac{1}{2} \binom{2^k}{2}} \prod_T M(T,T)^{-1/2} \prod_{\substack{(T,T')\\ T \neq T'}} M(T,T')^{-1/4}\right) \times \\ & \exp{\left\{ -n(d-1)H(\pi)\right\}} \times \prod_T M(T,T)^{-\frac{nd}{2}M(T,T)} \times \prod_{\substack{(T,T')\\ T \neq T'}} M(T,T')^{-\frac{nd}{2}M(T,T')}\\
&= O(\rm{poly}(n,d,M_{\min})) \exp{\left \{ n \, \left[\frac{d}{2}H(M) - (d-1)H(\pi)\right]\right\}}\,.
\qedhere\end{align*}
\end{proof}
The number of terms in the sum
$$
\E{Z(\rho)}=\sum_M \E{Z(\rho,M)}
$$
over edge profiles $M$ compatible with $\rho$ is bounded by a polynomial in $n$. Indeed, $M$ has $4^k$ non-negative entries of the form $m(T,T')/nd$ with the integers $m(T,T')$ satisfying $\sum_{(T,T')} m(T,T') = nd$. There are at most $(nd)^{4^k}$ such solutions. This allows us to conclude that $\E{Z(\rho)}$ is dominated by the largest exponential term, or in other words, the term with the largest value of $(d/2)H(M) - (d-1)H(\pi)$ optimized over $M$ that are compatible with $\rho$. We bound this optimum in the following.

Let $M = [ M(T,T')]_{\{T,T' \subset [k]\}}$ be an edge profile matrix with the property that $M$ is symmetric, the support of $M$ is contained in the
set $\{(T,T'): T \cap T' = \emptyset \}$ and that the marginal of $M$ along its row is a fixed probability distribution $\pi = (\pi(T); T \subset [k])$.
Define the weights $$w(T) = \sum_{T': T' \cap T = \emptyset} \pi(T').$$ Note that $w(\emptyset) = 1$, and
\begin{equation}
\label{wchosen}
\sum_{T: T \cap T' = \emptyset} \frac{\pi(T)}{w(T')} = 1.
\end{equation}

\begin{lem} \label{lem:maxentropy}
With a matrix $M$ and vectors $\pi,w$ as above
we have $$H(M) \leq 2H(\pi) + \sum_{S \subset [k]} \pi(T) \log(w(T))\,.$$
\end{lem}
\begin{proof}
Set $h(x) = - x\log(x)$ for $0 \leq x \leq 1 \; (0\log 0 = 0)$. Note that $h(x)$ is a smooth and strictly concave function on its domain.
We have that
\begin{eqnarray*}
H(M)
&=& \sum_{T'} \sum_{T: T \cap T' = \emptyset} \pi(T)
\frac{h(M(T,T'))}{\pi(T)}\quad
\\ &=& \sum_{T'} \sum_{T: T \cap T' = \emptyset} \pi(T) h\Big(\frac{M(T,T')}{\pi(T)}\Big) \quad + H(\pi)\,.
\end{eqnarray*}
For the second equality we used that $h(xy)=xh(y)+yh(x)$.

By Jensen's inequality applied to $h(x)$ and the identity \eqref{wchosen} we deduce that
\begin{eqnarray*}
\sum_{T: T \cap T' = \emptyset} \frac{\pi(T)}{w(T')} h\Big(\frac{M(T,T')}{\pi(T)}\Big) \leq h\Big(\sum_{T: T \cap T' = \emptyset} \frac{M(T,T')}{w(T')}\Big) = h\Big(\frac{\pi(T')}{w(T')}\Big)\,.\end{eqnarray*}
From this we conclude that
\[ H(M) \leq \sum_{T'} w(T') h(\frac{\pi(T')}{w(T')})  + H(\pi) = 2H(\pi) + \sum_{T'} \pi(T') \log(w(T'))\,.\qedhere\]
\end{proof}

Using Lemma \ref{lem:maxentropy} and Theorem \ref{thm:expectation} we conclude that for any density profile $\rho$
\begin{equation} \label{eqn:expectationbound} \E{Z(\rho)} \leq \rm{poly}(n,d) \times \exp{\left \{n \, \left [H(\pi) - \frac{d}{2}\hat{H}(\pi) \right ] \right \}}\end{equation}
where $\hat{H}(\pi) = \sum_{T} \pi(T) \log(w(T))$, and $\rm{poly}(n,d)$ is a polynomial in $n$ and $d$ of degree at most $4^k$.

For the purposes of our analysis we will be interested in density profiles $\rho$ such that $\rho(T) \in [\rho_{|T|} - \eps, \rho_{|T|}]$ with $\rho_i = \alpha_{i,d,p} \frac{\log d}{d}$. To this end let us fix $1 = \rho_0 \geq \rho_1 \geq \ldots \geq \rho_k$ with $\rho_i = \alpha_{i,d,p} \frac{\log d}{d}$.
Define the density profile $\rho$ by $\rho(T) = \rho_{|T|}$ for $T \subset [k]$. Let $\pi$ denote the probability distribution associated to
$\rho$ as given by (\ref{eqn:PIEpi}). For $T \neq \emptyset$ define the quantities $\beta(T)$ by $\pi(T) = \beta(T) \frac{\log d}{d}$.
Note that $\pi(\emptyset) = 1 - [\sum_{T \neq \emptyset} \beta(T)] \frac{\log d}{d}$. By setting $\alpha(T) = \alpha_{|T|,d,p}$
and using the relation between $\rho$ and $\pi$ from (\ref{eqn:PIEpi}) and (\ref{eqn:PIErho}) we conclude the following relation between
$\alpha$ and $\beta$:
\begin{eqnarray}
\label{eqn:PIEalpha} \alpha(T) &=& \sum_{T': T \subset T'} \beta(T) \\
\label{eqn:PIEbeta} \beta(T) &=& \sum_{T': T \subset T'} (-1)^{|T' \setminus T|} \alpha(T')\,.
\end{eqnarray}

Note that $\alpha_{1,d,p} = \alpha \leq 2$. Indeed, recall the result of Bollob\'{a}s \cite{Bol}
mentioned in the introduction: if $A_{n,d}$ is the event that all independent sets in $\G_{n,d}$ have
size at most $2 \frac{\log d}{d}\,n$ then $\pr{A_{n,d}} \to 1$ as $n \to \infty$ for every $d \geq 3$.
Recall from Section \ref{sec:projection} that $\alpha \frac{\log d}{d} = \lim_{n \to \infty} \E{|I_{\G_{n,d},1}|/n}$.
However, $\E{|I_{\G_{n,d},1}|} = \E{|I_{\G_{n,d},1}|; A_{n,d}} + \E{|I_{\G_{n,d},1}|; A_{n,c}^{c}} \leq 2 \frac{\log d}{d}\,n + n\pr{A_{n,d}^{c}}$.
Dividing through by $n$ and then taking limits in $n$, we deduce that $\alpha \leq 2$.

From the fact that $\alpha \leq 2$ we see that $0 \leq \alpha_{k,d,p} \leq \cdots \leq \alpha_{1,d,p} \leq 2$.
From (\ref{eqn:PIEbeta}) it follows that $\beta(T) \leq 2^{k+1}$ for all $T \subset [k]$.
In particular, this estimate is uniform in $d$ and $p$.

\begin{lem} \label{lem:asymptotics}
With $\pi$, $\alpha$ and $\beta$ as above we have that
$$H(\pi) - \frac{d}{2}\hat{H}(\pi) \leq \left [\sum_{i=1}^k (-1)^{i-1}\binom{k}{i} \alpha_{i,d,p}(2-\alpha_{i,d,p}) \right ] \frac{\log^2 d}{2d}
+ O_k(\frac{\log d}{d})$$
where the big $O$ term depends only on $k$.
\end{lem}

\begin{proof}
We need the asymptotic behaviour of $H(\pi) - \frac{d}{2}\hat{H}(\pi)$ where the entries of $\pi$ are on the scale of $(\log d)/d$.
By definition, 
$$
w(T) = 1 - \frac{\log d}{d}\sum_{T': T' \cap T \neq \emptyset} \beta(T').
$$ 
From Taylor expansion we observe that $-\log(1-x) \geq x$. Hence for $T \neq \emptyset$ we have
$$-\pi(T) \log(w(T)) \geq \beta(T) \left(\frac{\log d}{d}\right)^2 \sum_{T': T' \cap T \neq \emptyset} \beta(T')\,.$$
Since $w(\emptyset) = 1$ we have
$$\hat{H}(\pi) = \sum_{T \neq \emptyset} -\pi(T) \log(w(T)) \geq (\frac{\log d}{d})^2 \sum_{(T,T'): T \cap T' \neq \emptyset} \beta(T)\beta(T')].$$

To analyze $H(\pi)$ we consider the terms $h(\pi(\emptyset))$ and $h(\pi(T))$ with $T \neq \emptyset$ separately.
We note from Taylor expansion that $h(1-x) \leq x$ for $0 \leq x \leq 1$. Thus,

$$h(\pi(\emptyset)) = h\Big(1 - \frac{\log d}{d}\sum_{T \neq \emptyset} \beta(T) \Big) \leq
\frac{\log d}{d}\sum_{T \neq \emptyset} \beta(T)\,.$$

Since $\beta(T) \leq 2^{k+1}$ for $T \neq \emptyset$, we see that $h(\pi(\emptyset)) = O_k(\frac{\log d}{d})$.

On the other hand, for $T \neq \emptyset$ the quantity $h(\pi(T))$ equals 
$$h(\beta(T)\frac{\log d}{d}) = \beta(T)h(\frac{\log d}{d}) + h(\beta(T))\frac{\log d}{d} \leq \beta(T)\frac{\log^2 d}{d} + \frac{1}{e} \cdot \frac{\log d}{d}.$$
The inequality follows because $h(\frac{\log d}{d}) \leq \frac{\log^2 d}{d}$ and $h(x) \leq 1/e$ for all $x \geq 0$.

Therefore, $H(\pi) \leq \frac{\log^2 d}{d}\sum_{T \neq \emptyset} \beta(T) + O_k(\frac{\log d}{d})$.

From the above we conclude that
$$H(\pi) - \frac{d}{2}\hat{H}(\pi) \leq \frac{\log^2 d}{d} \Big [\sum_{T \neq \emptyset} \beta(T) - \frac{1}{2}\sum_{(T,T'): T \cap T' \neq \emptyset} \beta(T)\beta(T') \Big] + O_k(\frac{\log d}{d}).$$

Finally, it follows by inclusion-exclusion that 
$$\sum_{T \neq \emptyset} \beta(T) - \frac{1}{2}\sum_{(T,T'): T \cap T' \neq \emptyset} \beta(T)\beta(T') = \frac{1}{2} \sum_{i=1}^k (-1)^{i-1}\binom{k}{i} \alpha_{i,d,p}(2-\alpha_{i,d,p}).$$

The details are as follows. From the relations between $\alpha$ and $\beta$ in (\ref{eqn:PIEbeta}) and (\ref{eqn:PIEalpha}) it follows immediately that
$$\sum_{i=1}^k (-1)^{i-1}\binom{k}{i} \alpha_{i,d,p} = \sum_{T \neq \emptyset} \beta(T)$$
because both terms equal $(|\cup_{i=1}^k I_{G,i}|/n) \cdot \frac{d}{\log d}\,$.

Also, from these relations it follows that $\alpha(T)^2 = \sum_{(T_1, T_2): T \subset T_1 \cap T_2} \beta(T_1) \beta(T_2)$. Hence,
\begin{eqnarray*}
\sum_{i=1}^k (-1)^{i-1} \binom{k}{i} \alpha_{i,d,p}^2 &=& \sum_{T \neq \emptyset} (-1)^{|T|-1}\alpha(T)^2 \\
&=& \sum_{T \neq \emptyset} (-1)^{|T|-1} \sum_{(T_1, T_2): T \subset T_1 \cap T_2} \beta(T_1) \beta(T_2) \\
&=& \sum_{(T_1,T_2): T_1 \cap T_2 \neq \emptyset} \beta(T_1) \beta(T_2) \sum_{T: T \subset T_1 \cap T_2, T \neq \emptyset} (-1)^{|T|-1}\\
&=& \sum_{(T_1,T_2): T_1 \cap T_2 \neq \emptyset} \beta(T_1) \beta(T_2) \sum_{i=1}^{|T_1 \cap T_2|} (-1)^{i-1}\binom{|T_1 \cap T_2|}{i}\,.
\end{eqnarray*}

Now recall the binomial identity $\sum_{i=1}^t (-1)^{i-1} \binom{t}{i} = 1 - (1-1)^t = 1$ for any integer $t \geq 1$. This identity implies that
$$\sum_{i=1}^k (-1)^{i-1} \binom{k}{i} \alpha_{i,d,p}^2 = \sum_{(T,T'): T \cap T' \neq \emptyset} \beta(T) \beta(T')\,.$$
With this the proof of the final claim is complete.
\end{proof}

Let $E_{d,p}(\eps) = E(\alpha, \eps,n,d,p)$ be the event that $\G_{n,d}$ contains some $k$-tuple of independent sets
$(I_1, \ldots, I_k)$ whose density profile $\rho$ satisfies the property that
for every $T \subset [k]$,
\begin{equation} \label{eqn:Eevent}
\rho(T) \in \left [ \alpha_{|T|,d,p}\frac{\log d}{d} - \eps, \alpha_{|T|,d,p} \frac{\log d}{d} + \eps \right].
\end{equation}

We can bound $\pr{E_{d,p}(\eps)}$ from above via (\ref{eqn:expectationbound}) and Lemma \ref{lem:asymptotics}.
Define the density profile $\rho_{\alpha}$ by $\rho_{\alpha}(T) = \alpha_{|T|,d,p} \frac{\log d}{d}$. Let
$\pi_{\rho_{\alpha}}$ be the corresponding probability vector obtained from (\ref{eqn:PIEpi}).
For any admissible density profile $\rho$ for the occurrence of the event
$E_{d,p}(\eps)$, the corresponding $\pi_{\rho}$ satisfies $|\pi_{\rho}(T) - \pi_{\rho_{\alpha}}| = O_k(\eps)$.
We employ Lemma \ref{lem:asymptotics} for $\pi_{\rho}$. We get that
$$H(\pi_{\rho}) - \frac{d}{2}\hat{H}(\pi_{\rho})= H(\pi_{\rho_{\alpha}}) - \frac{d}{2}\hat{H}(\pi_{\rho_{\alpha}}) + \rm{err}_{d,k}(\eps)\,.$$

The error term $\rm{err}_{d,k}$ is such that $\rm{err}_{d,k}(\eps) \to 0$ as $\eps \to 0$, and this holds uniformly in $p \in [0,1]$.
This follows from the fact that $\pi$ is obtained from $\rho$ by a smooth transformation (see (\ref{eqn:PIEpi})),
and that $H$ and $\hat{H}$ are smooth functions. The reason $\rm{err}_{d,k}(\eps)$
tends to 0 uniformly in $p$ is because it depends on the $\alpha_{i,d,p}$ smoothly and only through their absolute values.
However, the $\alpha_{i,d,p}$ are all bounded as $0 \leq \alpha_{k,d,p} \leq \cdots \leq \alpha_{1,d,p} = \alpha \leq 2$.
A careful analysis will actually show that $\rm{err}_{d,k}(\eps) = O_k(\frac{\log^2 d}{d} \,\eps)$.

From Lemma \ref{lem:asymptotics} applied to $\rho_{\alpha}$ it follows that for any admissible $\rho$ for the occurrence of the event
$E_{d,p}(\eps)$,
\begin{equation} \label{eqn:rhobound}
H(\pi_{\rho}) - \frac{d}{2} \hat{H}(\pi_{\rho}) \leq \frac{\log^2 d}{d} \sum_{i=1}^k (-1)^{i-1} \binom{k}{i} \alpha_{i,d,p}(2-\alpha_{i,d,p}) + O_k(\frac{\log d}{d}) + \rm{err}_{d,k}(\eps)\,.\end{equation}

Now note that the number of density profiles $\rho$ that is admissible for the event $E_{d,p}(\eps)$ is at most $O(n^{2^k})$ where the big O constant is uniformly bounded in $n$ because quantities of the form $\alpha_{i,d,p}\frac{\log d}{d}$ are all of constant order in $n$. Taking an union bound over all such admissible $\rho$, using the first moment method and employing the bounds in (\ref{eqn:expectationbound}) and (\ref{eqn:rhobound}), we conclude that $\pr{E_{d,p}(\eps)}$ is bounded above by a polynomial term $\rm{poly}(n,d)$ times the exponential term
\begin{equation} \label{eqn:Eupbound}
\exp {\left \{ n\left [\frac{\log^2 d}{d} \Big(\sum_{i=1}^k (-1)^{i-1} \binom{k}{i} \alpha_{i,d,p}(2-\alpha_{i,d,p})\Big) + O_k(\frac{\log d}{d}) + \rm{err}_{d,k}(\eps) \right]\right\}} \,.
\end{equation}

\subsection{Concentration of the density profile of $(I_{\G,1}, \ldots, I_{\G,k})$ about its mean} \label{sec:concentration}

Having derived an upper bound to $\pr{E_{d,p}(\eps)}$ we need a lower bound on this probability
that violates the upper bound and provides a contradiction. We now show that $\pr{E_{d,p}(\eps)} \to 1$
as $n \to \infty$ via concentration inequalities.

Recall the terminology of Section  \ref{sec:binomineq}. The factor of i.i.d.~independent set $I_d$ is used to construct independent sets
$I_{G,1}, \ldots, I_{G,k}$ of $G \in G_{n,d}$ using random labellings $X_0, \ldots, X_k$ of $G$ and a random subset $S \subset V(G)$
resulting from Bernoulli percolation on $G$. Set $Y(v) = (X_0(v), \ldots, X_k(v), \mathbf{1}_{v \in S})$ for $v \in V(G)$. 
Then $Y = (Y(v); v \in V(G))$ is an i.i.d.~process on $G$ and $I_{G,1}, \ldots, I_{G,k}$ is a function of $Y$. Also, for any $T \subset [k]$
the expected density $\E{|\cap_{i \in T} I_{G,i}|/n} \in [\alpha_{i,d,p} \frac{\log d}{d}(1- \frac{B(G)}{n}), \alpha_{i,d,p} \frac{\log d}{d}]$,
where $B(G)$ is the number of vertices in $G$ whose $(r+1)$-neighbourhood is not a tree.

\begin{lem} \label{lem:concentration}
For any $G \in G_{n,d}$ the independent sets $I_{G,1}, \ldots, I_{G,k}$ satisfy the following with $C_{r,d} = O(r^2 d^{2r})$:
\begin{equation}
\pr{\max_{T \subset [k]} \left | \frac{|\cap_{i \in T} I_{G,i}| - \E{|\cap_{i \in T} I_{G,i}|}}{n} \right | > \eps} \leq 2^{k+1} e^{-\frac{\eps^2 n}{C_{r,d}}}\,.
\end{equation}
\end{lem}

\begin{proof} For each $T \subset [k]$ the set $\cap_{i \in T} I_{G,i}$ is a function of $y = (y(v); v \in V(G))$,
where each $y(v) \in [0,1]^{k+1} \times \{0,1\}$ (the set of values of the random variable $Y(v)$).
Modifying some entry $y(v)$ to $y'(v)$ can switch the state of inclusion of a vertex $u$ within $\cap_{i \in T} I_{G,i}$ only if
$u$ is in $N_{G}(r,v)$, where $r$ is the radius of the factor associated to $I_d$. Therefore, such a modification to $y$
can cause the size of $\cap_{i \in T} I_{G,i}$ to change by at most $|N_{G}(r,v)| = O(rd^r)$ since $G$ is $d$-regular.
Since the random input $Y$ is an i.i.d.~process it follows from the Hoeffding--Azuma inequality \cite[Theorem 1.20]{Bolbook} that
$$\pr{\Big | |\cap_{i \in T} I_{G,i}| - \E{|\cap_{i \in T} I_{G,i}|} \Big | > x} \leq 2e^{\frac{x^2}{2nC_{r,d}}}\,.$$
The lemma follows by taking an union bound over $T \subset [k]$ and replacing $x$ by $n\eps$.
\end{proof}

Recall that for the random graph $\G_{n,d}$ we have
$$\alpha_{|T|,d,p} \,\frac{\log d}{d}(1 - \frac{\E{B(\G_{n,d})}}{n}) \leq  \E{\frac{|\bigcap_{i \in T} I_{\G_{n,d},i}|}{n}} \leq \alpha_{|T|,d,p} \,\frac{\log d}{d}\,.$$
As we mentioned in Section \ref{sec:projection}, $\E{B(\G_{n,d})} = O_{d,r}(1)$ in $n$.

We may find an $n_d$ such that $\alpha_{k,d,p} \frac{\log d}{d} \, \E{B(\G_{n,d})} \leq (\eps/2) n$ for $n \geq n_d$.
This ensures that
$$\alpha_{|T|,d,p} \, \frac{\log d}{d} \, \left(1 - \frac{\E{B(\G_{n,d})}}{n} \right) \geq \alpha_{|T|,d,p} \, \frac{\log d}{d} - \frac{\eps}{2}
\quad \text{for all} \quad T\subset [k].$$
If $\big | |\cap_{i \in T} I_{\G_{n,d},i}| - \E{|\cap_{i \in T} I_{\G_{n,d},i}|} \big | \leq (\eps/2)n$ then the
independent sets $I_{\G_{n,d},i}$ satisfy, for any $T \subset [k]$,

$$\alpha_{|T|,d,p} \frac{\log d}{d} - \eps \leq \frac{|\cap_{i \in T} I_{\G_{n,d},i}|}{n} \leq \alpha_{|T|,d,p} \frac{\log d}{d} + \eps.$$

Therefore, the event $E_{d,p}$ occurs for $n \geq n_d$ (see the definition of $E_{d,p}$ in (\ref{eqn:Eevent})).
From Lemma \ref{lem:concentration} we conclude that
for $n \geq n_d$,
\begin{equation} \label{eqn:Ehighprob}
\pr{E_{d,p}(\eps)} \geq 1 - 2^{k+1} e^{-\frac{\eps^2}{C_{r,d}}\,n} \longrightarrow 1 \quad \text{as}\; n \to \infty\,.
\end{equation}

\subsubsection{Conclusion of the proof of Theorem \ref{thm:binomineq}} \label{sec:regulargraphconclusion}

Suppose for some $\delta > 0$ we have that
$$\liminf_{d \to \infty} \, \inf_{p \in [0,1]}\, \sum_{i=1}^k (-1)^{i-1} \binom{k}{i} \alpha_{i,d,p}(2-\alpha_{i,d,p}) = -\delta.$$

For each $d$ we pick a $p' = p'(d)$ such that
$$\inf_{p \in [0,1]}\, \sum_{i=1}^k (-1)^{i-1} \binom{k}{i} \alpha_{i,d,p}(2-\alpha_{i,d,p}) \geq
 \sum_{i=1}^k (-1)^{i-1} \binom{k}{i} \alpha_{i,d,p'}(2-\alpha_{i,d,p'}) - \frac{\delta}{2}.$$
 
We deduce that,
\begin{equation} \label{eqn:liminfneg}
\liminf_{d \to \infty} \; \sum_{i=1}^k (-1)^{i-1} \binom{k}{i} \alpha_{i,d,p'}(2-\alpha_{i,d,p'}) \leq - \frac{\delta}{2}.
\end{equation}

Recall the upper bound for $\pr{E_{d,p}(\eps)}$ in (\ref{eqn:Eupbound}), which we now consider for $p = p'(d)$.
The error term $\rm{err}_{d,k}(\eps) \to 0$ as $\eps \to 0$, uniformly in $p$. So for each $d$ pick an $\eps_d$ such that
$\rm{err}_{d,k}(\eps_d) \leq (\delta/4) \frac{\log^2 d}{d}$. The inequalities (\ref{eqn:Eupbound}) and (\ref{eqn:liminfneg})
imply that there exists a subsequence $d_i \to \infty$ such that
$$\pr{E_{d_i, p'(d_i)}(\eps_d)} \leq \mathrm{poly}(n,d_i) \, e^{n \left [-\frac{\delta \log^2 d_i}{4d_i} \,+\, O_k \big (\frac{\log d_i}{d_i} \big) \right ]}
\longrightarrow 0 \;\;\text{as}\; n \to \infty$$
for all sufficiently large values of $d_i$.
However, we have already concluded from (\ref{eqn:Ehighprob}) that $\pr{E_{d_i, p'(d_i)}(\eps_d)} \to 1$ as $n \to \infty$ for all such $d_i$.
This provides a contradiction and completes the proof.

\section{Local algorithms for independent sets in Erd\H{o}s-R\'{e}nyi graphs} \label{sec:pgwtrees}

Local algorithms on sparse Erd\H{o}s-R\'{e}nyi graphs are projections of factor of i.i.d.~processes on Poisson-Galton-Watson (PGW) trees.
We will define the appropriate notion of factor of i.i.d.~independent sets in PGW trees and prove the same asymptotic upper and lower bounds as for regular trees. Recall that a PGW tree with average degree $\lambda$, which we denote $\pgw$, is a random tree resulting from a Galton-Watson branching process with a $\poi(\lambda)$ offspring distribution. Before we can define the notion of factor of i.i.d.~independent sets in PGW trees we will need some notation.

Let $\Lambda_r$ denote the collection of all triples $(H,v,x)$ where (1) $(H,v)$ is a finite, connected, rooted graph with root $v$, (2) for all vertices $u \in V(H)$ we have $dist(u,v) \leq r$ where $dist$ denotes the graph distance, and (3) $x \in [0,1]^{V(H)}$ is a labelling of $H$. $\Lambda_r$ has a natural $\sigma$-algebra, $\Sigma_r$, generated by sets of the form $(H,v)\times B$ where $(H,v)$ satisfies properties (1) and (2) above and $B$ is a Borel measurable subset of $[0,1]^{V(H)}$. We consider two rooted graphs to be isomorphic if there exists a graph isomorphism between them that maps one root to the other. Given an isomorphism $\phi : (H,v) \to (H',v')$, any labelling $x$ of $(H,v)$ induces a labelling $\phi \cdot x$ of $(H',v')$ by defining $\phi \cdot x(i) = x(\phi^{-1}(i))$, and vice-versa. A function $f : \Lambda_r \to \{0,1\}$ is a factor if it is $\Sigma_r$ measurable and $f(H,v,x) = f(\phi(H), \phi(v), \phi \cdot x)$ for all isomorphisms $\phi$ of $H$, and all $H$.

For $0 \leq r < \infty$, let $f: \Lambda_r \to \{0,1\}$ be a factor. Consider a $\pgw$ tree with a random labelling $X$. Let $N_r(\pgw,v)$ denote the $r$-neighbourhood of a vertex $v$ in $\pgw$ and let $X(\pgw,v,r)$ be the restriction of $X$ to $N_r(\pgw,v)$. Define a subset $I$ of the vertices of $\pgw$ by setting
$$I(v) = f(N_r(\pgw,v),v, X(\pgw,v,r))\,.$$
We say that $I$ is a factor of i.i.d.~independent set in $\pgw$ if $I$ is an independent set in this tree with probability 1 (w.r.t.~the random labelled tree $(\pgw,X)$).

The distribution of the random variable $I(v)$ does not depend on the choice of the vertex $v$. This is because in a PGW tree the distribution of the neighbourhoods $N_r(\pgw,v)$ does not depend on the choice of $v$. So let $\mathrm{PGW}(\lm,r)$ denote the tree following the common distribution of these $r$-neighbourhoods, rooted at a vertex $\circ$, and let $X$ be a random labelling. The density of the subset $I$ is defined to be the expectation
$$density(I) = \E{f(\mathrm{PGW}(\lm,r),\circ,X)}\,.$$

Define the quantity $\alpha(\lm)$ by
$$\alpha(\lm) \frac{\log \lm}{\lm} =
\sup_{0 \leq r < \infty} \left \{ density (I) : I \;\text{an independent set in}\; \pgw \;\text{with factor}\; f : \Lambda_r \to \{0,1\} \right \}\,.$$

\begin{thm} \label{thm:pgwthm}
The limit $\lim_{\lm \to \infty} \alpha(\lm) = 1$.
\end{thm}

In Section \ref{sec:pgwupbound} we prove that $\limsup_{\lm \to \infty} \alpha(\lm) \leq 1$, and in Section \ref{sec:pgwlwbound} that $\liminf_{\lm \to \infty} \alpha(\lm) \geq 1$. The proof of the upper bound will employ the strategy used for regular trees in Section \ref{sec:regulartrees}. We will highlight the key differences but be brief with parts of the argument that are analogous to the case for regular trees.

\subsection{Upper bound on density of factor of i.i.d.~independent sets in PGW trees} \label{sec:pgwupbound}

Recall the Erd\H{o}s-R\'{e}nyi graph $ER(n, p)$ is a random graph on the vertex set $[n]$ where every pair of vertices $\{u,v\}$ is independently included with probability $p$. Our interest lies with the random graphs $ER(n, \lambda/n)$ where $ \lambda > 0$ is fixed. Throughout this section let $G_n$ denote a random graph that is distributed according to the $ER(n,\lm/n)$ model. It is well known (see \cite{Bolbook} chapter 4) that the sequence of random graphs $G_n$ converges in the \emph{local weak limit} to the tree $\pgw$. This means that for every fixed $r \geq 0$, if $\circ_n \in [n]$ is chosen uniformly at random then for any finite rooted graph $(H,\circ)$ the probability $\pr{N_r(G_n,\circ_n) \cong (H,\circ)} \to \pr{\mathrm{PGW}(\lm,r) \cong (H,\circ)}$ as $n \to \infty$.

Consequently, using the same technique that was used for regular trees, a factor of i.i.d.~independent set $I$ of $\pgw$ with factor $f: \Lambda_r \to \{0,1\}$ yields a factor of i.i.d.~independent set $I_n$ of $G_n$ in the following sense. These exists a factor $f_n : \Lambda_{r+1} \to \{0,1\}$ such that if $X$ is a random labelling of $G_n$ then $I_n(v) = f_n(N_{r+1}(G_n,v),v, X(G_n,v,r+1))$. Furthermore, $\E{|I_n|/n} \to density(I)$ as $n \to \infty$.

To prove that $\limsup \alpha(\lm) \leq 1$ we assume to the contrary. Then we can find $\alpha > 1$ and a subsequence of $\lm \to \infty$ such that for each $\lm$ there exists a factor of i.i.d.~independent set $I_{n,\lm}$ of $G_n$ with factor $f_{n,\lm} : \Lambda_{r_{\lm}} \to \{0,1\}$, and $E{|I_{n,\lm}|/n} \geq \alpha \frac{\log \lm}{\lm}$ for all sufficiently large $n$. We can assume w.l.o.g.~that these statements hold for all $\lm$ and $n$. By setting $E{|I_{n,\lm}|/n} = \alpha_{1,n,\lm} \frac{\log \lm}{\lm}$ we have that $\alpha_{1,n,\lm} \geq \alpha > 1$.

\subsubsection{A coupling of local algorithms on Erd\H{o}s-R\'{e}nyi graphs} \label{sec:ERcoupling}

For $0 \leq p \leq 1$ let $S = S_{n,p}$ be a random subset of $V(G_n)$ chosen by doing a Bernoulli percolation with density $p$. Let $G'_n = G'_n(G_n,S)$ be the random graph that is obtained from $G_n$ by independently resampling the edge connections between each pair of vertices $\{u,v\} \subset S$ with inclusion probability $\lm/n$. In other words, $G'_n$ retains all edges of $G_n$ that do not connect $S$ to itself, and all possible edge connections between vertices within $S$ are resampled according to the Erd\H{o}s-R\'{e}nyi model. Note that $G'_n$ is also distributed according to $ER(n, \lm/n)$; if $p=0$ then $G'_n = G_n$, and if $p=1$ then $G'_n$ is independent of $G_n$.

Now fix $G_n$ and $S$ as above and let $X$ be a random labelling of $G_n$. Let $X^1, X^2, \ldots$ be new, independent random labellings and define labellings $Y^k$, correlated with $X$, by $Y^k(v) = Y^k(v)$ if $v \in S$, and $Y^k(v) = X(v)$ otherwise. Generate $G^1,G^2,\ldots$ from $G_n$ and $S$
by using the recipe for $G'_n$, but rewire the edges for each $G^i$ independently. In other words, the induced subgraphs $G^1[S], G^2[S],\ldots$ are independent. Now consider independent sets $I^1,I^2,\ldots$ by letting $I^k$ be generated by the factor $f_{n,\lm}$ with input graph $G^k$ and labelling $Y^k$. Thus, $I^k$ is a factor of i..i.d. independent set of $G^k$. Since all these graphs have a common vertex set, namely $[n]$, we can consider intersections of the $I^k$. Note that for any finite subset $T$ the expected intersection density $\E{|\cap_{t \in T} I^t|/n}$ depends only of $|T|$ due to exchangeability of the $I^k$. Define the parameters $\alpha_{k,n,\lm}$ by
\begin{equation} \label{eqn:PGWintdensity}
\alpha_{k,n,\lm} \frac{\log \lm}{\lm} = \E{\frac{|\cap_{t=1}^k I^t|}{n}}.
\end{equation}

\begin{thm} \label{thm:pgwbinomeq}
The following inequality holds for each $k \geq 1$
\begin{equation} \label{pgwbinomeq}
\liminf_{\lm \to \infty} \, \liminf_{n \to \infty} \, \inf_{p \in [0,1]} \, \sum_{i=1}^k (-1)^{i-1}\binom{k}{i} \alpha_{i,n,\lm,p}(2-\alpha_{i,n,\lm,p}) \geq 0\,.
\end{equation}
\end{thm}

Now we establish the upper bound by using Theorem \ref{thm:pgwbinomeq}. We define stability variables $Q_{n,\lm,p}$ in a manner analogous
to what we did for regular graphs. First, let $\circ \in [n]$ be a uniform random vertex. We restrict our probability space to the support of
$f_{n,\lm}(N_{r_{\lm}}(G_n,\circ),\circ,X(G_n,\circ,r_{\lm}))$, considered as a subset of the original probability space
determined by $\circ$, the random labellings $X, X^1,\ldots$, the random subset $S$ and the independent trials that determine the graphs
$G_n, G^1, G^2,\ldots$. Let $\mathbb{E}^{*}$ be the expectation operator $\mathbb{E}$ restricted to the new space:
$$ \mathbb{E}^{*}\left [ U \right ] =
\frac{\E{U \, f_{n,\lm}(N_{r_{\lm}}(G_n,\circ),\circ,X(G_n,\circ,r_{\lm}))}}{\E{f_{n,\lm}(N_{r_{\lm}}(G_n,\circ),\circ,X(G_n,\circ,r_{\lm}))}}\,.$$
Notice that we define the new probability space on finite graphs instead of on the infinite limiting graph as we did previously for regular graphs.
This coupling takes into account the randomness in the local structure of the underlying Erd\H{o}s-R\'{e}nyi graphs, which is not an issue for regular graphs.

Define the stability $Q_{n,\lm,p} = Q_{n,\lm}(G_n,\circ,S,X)$ on the new probability space by
$$ Q_{n,\lm,p} = \Est{f_{n,\lm}(N_{r_{\lm}}(G^1,\circ),\circ, Y^1(G^1,\circ, r_{\lm}) \mid G_n,\circ,S,X})\,.$$
One can check, as before, that for every $k \geq 1$ the moment $\Est{Q^{k-1}_{n,\lm}} = \frac{\alpha_{k,n,\lm,p}}{\alpha_{1,n,\lm}}$.

We now show that expectations involving $Q_{n,\lm,p}$ are continuous with respect to $p$, and that $Q_{n,\lm,p}$ has the right values at the endpoints $p=0$ and $p=1$. Observe that $Q_{n,\lm,p} \in [0,1]$. If $g$ is a bounded measurable function on $[0,1]$ then $\Est{g(Q_{n,\lm,p})}$ is Lipschitz in $p$.

Indeed, let $p_1 \leq p_2$. We couple the labelled graphs $(G^1(S_{p_1}),\circ,Y^1_{p_1})$ and $(G^1(S_{p_2}), \circ, Y^1_{p_2})$
given $(G_n,\circ,X)$ through the percolation subsets. Let $Z$ be a random labelling of $[n]$,
and let $\tau_{\{u,v\}}$ for $\{u,v\} \subset [n]$ be independent Bernoulli trials of expectation $\lm/n$. Set $S_{p_1} = \{ v: Z(v) \leq p_1\}$ and $S_{p_2} = \{ v : Z(v) \leq p_2\}$. The resampled edges of $G^1(S_{p_1})$ (resp. $G^1(S_{p_2})$) are determined according to the $\tau_{\{u,v\}}$ for $u,v \in S_{p_1}$ (resp. for $u,v \in S_{p_2}$). Similarly, the labelling $Y^1_{p_1}$ (resp. $Y^1_{p_2}$) agrees with $X^1$ on $S_{p_1}$ (resp. $S_{p_2}$) and agrees with $X$ otherwise. With this coupling we have that (ignoring some formalities with the notation)
\begin{align*}
&\Est{g(Q_{n,\lm,p_1})} - \Est{g(Q_{n,\lm,p_2})}  \; = \\
&\Est{g\left(\Est{f(G^1(S_{p_1}), \circ, Y^1_{p_1}) \mid G_n,\circ, X,S_{p_1}}\right)
- g\left(\Est{f(G^1(S_{p_2}), \circ, Y^1_{p_2}) \mid G_n,\circ, X,S_{p_2}}\right)}.
\end{align*}

Observe that if $Z(v) \notin (p_1,p_2)$ for every $v \in [n]$ then $S_{p_1} = S_{p_2}$, and hence, $(G^1(S_{p_1}), Y^1_{p_1}) = (G^1(S_{p_2}), Y^1_{p_2})$.
On this event, $\Est{f(G^1(S_{p_1}), \circ, Y^1_{p_1}) \mid G_n,\circ, X,S_{p_1}}$ equals \newline
$\Est{f(G^1(S_{p_2}), \circ, Y^1_{p_2}) \mid G_n,\circ, X,S_{p_2}}$.
So the difference of the two expectations above is zero on this event.
By an union bound, the probability that $Z(v) \in (p_1,p_2)$ for some vertex $v$ is at most $n|p_1-p_2|$. Therefore, it follows
from the triangle inequality that
$\left |\, \Est{g(Q_{n,\lm,p_1})} - \Est{g(Q_{n,\lm},p_2)} \,\right | \leq (2||g||_{\infty}n)\,|p_1-p_2|.$

The endpoint values of $Q_{n,\lm,p}$ are the same as before. When $p=0$ the resampled graph $G^1$ equals $G_n$, and the labelling $Y^1 = X$ due to $S$ being empty. Consequently $Q_{n,\lm,0} \equiv 1$ on the restricted probability space. On the other hand, if $p=1$ then $(G^1,Y^1)$ is independent of $(G_n,X)$ and the conditioning has no effect due to $S$ being the entire vertex set. Note that the common root $\circ$ does not affect the calculation because the 
distribution of $N_r(\mathrm{ER}(n,\lm/n),v)$ does not depend on $v$. We thus have $Q_{n,\lm,1} \equiv \alpha_{1,n,\lm} \frac{\log \lm}{\lm}$.

With these observations we can now proceed with the proof exactly the same way as before.
We skip the remainder of the argument for brevity and prove Theorem \ref{thm:pgwbinomeq} in the following.

\subsubsection{Proof of Theorem \ref{thm:pgwbinomeq}}

We will show that the existence of the factor of i.i.d.~independent sets $I^i$ on the graph $G^i$
implies that with high probability each graph $G^i$ contains a subset $S^i$ such that $S^i$ is an independent set in $G^i$,
and the empirical intersection densities of the $S^1, \ldots, S^k$ are close to the quantities $\alpha_{k,n,\lm} \frac{\log \lm}{\lm}$
Then we will bound the probability of observing such a $k$-tuple of independent sets, and prove that this probability is vanishingly
small unless Theorem (\ref{thm:pgwbinomeq}) holds.

Consider subsets $S^i \subset V(G^i)$. The density profile of the $k$-tuple $(S^1, \ldots, S^k)$ is
the vector $\rho = (\rho(T); T \subset [k])$ where $\rho = (\rho(T); T \subset [k])$ defined as
$\rho(T) = \frac{|\cap_{i \in T} S^i|}{n}$. Consider the independent sets $I^i$ of $G^i$ defined in Section \ref{sec:ERcoupling}.
They satisfy $\E{|\cap_{j=1}^i I^j |/n} = \alpha_{i,n,\lm,p} \frac{\log \lm}{\lm}$ for every $1 \leq i \leq k$.

Fix $0 < \eps < 1$. Let $A(\eps,p)$ be the following event. For each $1 \leq i \leq k$, $G^i$ contains
an independent set $S^i$ such that the density profile of $(S^1, \ldots, S^k)$ satisfies the following for all $T \subset [k]$:
$$\rho(T) \in [(1-\eps) \alpha_{|T|,n,\lm,p}\frac{\log \lm}{\lm}, (1+\eps) \alpha_{|T|, n, \lm,p}\frac{\log \lm}{\lm}].$$

We show that $\pr{A(\eps,p)} \to 1$ as $n \to \infty$. This follows if we show that
$$\pr{\max_{T \subset [k]} \left |\frac{|\cap_{i \in T} I^i|}{n} - \E{\frac{|\cap_{i \in T} I^i|}{n}} \right | > \eps} \to 0.$$
Indeed, given a realization of the graphs $G^1, \ldots, G^k$ and random labellings $Y^1, \ldots, Y^k$,
we take $S^i = I^i$ on $G^i$. If $|\frac{|\cap_{i \in T} I^i|}{n} - \E{\frac{|\cap_{i \in T} I^i|}{n}} | \leq \eps$
for every $T \subset [k]$ then the conditions for $A(\eps,p)$ to occur are satisfied.

\begin{lem} \label{lem:pgwconcentration}
With $G^1, \ldots, G^k$ as defined and corresponding independent sets $I^1, \ldots, I^k$ as defined via the factor
$f_{n,\lm}$, one has that for all $\eps > 0$, as $n \to \infty$,
$$\pr{\max_{T \subset [k]} \left |\frac{|\cap_{i \in T} I^i|}{n} - \E{\frac{|\cap_{i \in T} I^i|}{n}} \right | > \eps} \to 0\,.$$
\end{lem}

\begin{proof}
We show that $\E{\left | |\cap_{i \in T} I^i| - \E{|\cap_{i \in T} I^i|} \right |^2} = o(n^2)$ where the little o term may depend on $\lm, r_{\lm},$ and $k$. The statement of the lemma then follows from Chebyshev's inequality and an union bound over $T \subset [k]$. We write

$$ |\cap_{i \in T} I^i| - \E{|\cap_{i \in T} I^i|} = \sum_{v=1}^n \mathbf{1}\{v \in \cap_{i \in T} I^i\} - \pr{v \in \cap_{i \in T} I^i}\,.$$

Now, $| \mathbf{1}\{v \in \cap_{i \in T} I^i\} - \pr{v \in \cap_{i \in T} I^i}| \leq 2$ and so $\E{(\mathbf{1}\{v \in \cap_{i \in T} I^i\} - \pr{v \in \cap_{i \in T} I^i})^2} \leq 4$. Also, for two vertices $u$ and $v$ if the graph distance $dist_{G^i}(u,v) > 2r_{\lm}$ then the events $\{u \in I^i\}$ and $\{v \in I^i\}$ are independent with respect to the random labelling of $G^i$ because the factor $f_{n,\lm}$ makes decisions based on the labels along the $r_{\lm}$-neighbourhood of a vertex. Consequently,
$\E{(\mathbf{1}\{u \in \cap_{i \in T} I^i\} - \pr{u \in \cap_{i \in T} I^i})\cdot(\mathbf{1}\{v \in \cap_{i \in T} I^i\} - \pr{v \in \cap_{i \in T} I^i})}$ is at most $4 \pr{dist_{G^i}(u,v) \leq 2r_{\lm} \; \text{for some}\; i}.$

These two observations imply that
$$\E{\left | |\cap_{i \in T} I^i| - \E{|\cap_{i \in T} I^i|} \right |^2} \leq 4n + 4 \E{\# \left \{ (u,v): dist_{G^i}(u,v) \leq 2r_{\lm} \; \text{for some}\; i \right \}}.$$

Using the fact that the random graphs $ER(n,\lm/n)$ converge locally to $\pgw$, it is a standard exercise to show that the expected number of pairs $(u,v)$ in $ER(n,\lm/n)$ that satisfy $dist(u,v) > R$ is $o(n^2)$ (the little o term depends on $\lm$ and $R$). From this observation and a union bound over $i$ we deduce that $\E{\# \{(u,v): dist_{G^i}(u,v) \leq 2r_{\lm} \; \text{for some}\; i\}} = o(n^2)$. This proves the estimate for the squared expectation and completes the proof.
\end{proof}

In order to bound $\pr{A(\eps,p)}$ from above we need a procedure to sample the graphs $G^1, \ldots, G^k$.

\paragraph{\textbf{Sampling the graphs $(G^1, \ldots, G^k)$}} Let $\tau_{i,u,v}$ for $1 \leq i \leq k$ and $\{u,v\} \subset [n]$ be the indicator of the event that the edge $\{u,v\}$ belongs to $G^i$. Then the random vectors $(\tau_{i,u,v}; 1 \leq i \leq k)$ are independent of each other as $\{u,v\}$ varies. Let $S \subset [n]$ be a random subset chosen by a Bernoulli percolation with density $p$. If both $u,v \in S$ then $(\tau_{i,u,v}; 1 \leq i \leq k)$ are independent Bernoulli trials of expectation $\lm/n$ for each $i$. Otherwise, $(\tau_{i,u,v}; 1 \leq i \leq k)$ satisfies $\tau_{1,u,v} = \cdots = \tau_{k,u,v}$. In the latter case all $k$ of these indicators take the value 1 with probability $\lm/n$ or they are all zero with the complementary probability.

The sampling procedure above will allow us to compute expectations involving independent sets in the $G^i$. Let $I^i \subset [n]$ be an independent set of $G^i$. Defining $\rho(T) = |\cap_{t \in T} I^t|/n$ for $T \subset [k]$, the density profile associated to these $k$ independent sets is $\rho = (\rho(T); T \subset [k])$. The density profile $\rho$ determines a probability distribution $\pi = (\pi(T); T \subset [k])$ by equation (\ref{eqn:PIEpi}). Let $Z(\rho)$ be the number of $k$-tuple of subsets $(I^1,\ldots,I^k)$ of $[n]$ such that they have density profile $\rho$ and $I^i$ is an independent set of $G^i$.

\begin{lem} \label{lem:pgwexpectation}
The expectation of $Z(\rho)$ satisfies
$$\E{Z(\rho)} \leq \binom{n}{\pi(T)n\,; T \subset [k]} \left (1-\frac{\lm}{n}\right )^{\sum_{T \neq \emptyset} \binom{\pi(T)n}{2} + \frac{1}{2} \sum_{(T,T'): T \neq T', T\cap T' \neq \emptyset} \pi(T)\pi(T')n^2} \,.$$
\end{lem}

\begin{proof}
Fix a $k$-tuple $(I^1,\ldots,I^k)$ with each $I^i \subset [n]$ such that density profile of the $k$-tuple is $\rho$.
Given $v \in [n]$, let $T_v = \{ i \in [k]: v \in I^i\}$. Let $E_{\{u,v\}}$ be the event that the edge $\{u,v\}$ is absent is all $G^i$ for which $i \in T_u \cap T_v$, that is, $E_{\{u,v\}} = \{\tau_{i,u,v} = 0$ for all $i \in T_u \cap T_v\}$. The subsets $I^1, \ldots, I^k$ have the property that $I^i$ is an independent set of $G^i$ if and only if the events $E_{\{u,v\}}$ occur for all pairs $\{u,v\}$.

From the sampling procedure for the graphs $G^1, \ldots, G^k$, we note that the events $E_{\{u,v\}}$ are independent. Conditioning on the random subset $S$ and using the sampling procedure we conclude that
$$\E{Z(\rho)|S} = \sum_{\substack{k-tuples \, (I^1,\ldots,I^k) \\ \text{with density profile}\; \rho}}
\prod_{\{u,v\} \subset S} (1-\frac{\lm}{n})^{|T_u \cap T_v|} \,\prod_{\text{all other}\; \{u,v\}} (1-\frac{\lm}{n})^{\mathbf{1}_{\{T_u \cap T_v \neq \emptyset \}}}\,.$$

Observe that $(1-\frac{\lm}{n})^{|T_u \cap T_v|} \leq (1-\frac{\lm}{n})^{\mathbf{1}_{\{T_u \cap T_v \neq \emptyset \}}}$. Therefore, no matter the outcome of $S$ we have that
$$\E{Z(\rho)|S} \leq  \sum_{\substack{k-tuples \, (I^1,\ldots,I^k) \\ \text{with density profile}\; \rho}}
\prod_{\{u,v\}} (1-\frac{\lm}{n})^{\mathbf{1}_{\{T_u \cap T_v \neq \emptyset \}}}\,.$$

This implies the same inequality for the unconditional expectation $\E{Z(\rho)}$. The key observation is that
\begin{equation} \label{eqn:intersectioncount}
\sum_{\{u,v\}} \mathbf{1}_{\{T_u \cap T_v \neq \emptyset \}} = \sum_{T \neq \emptyset} \binom{\pi(T)n}{2} + \frac{1}{2} \sum_{(T,T'): T \neq T', T\cap T' \neq \emptyset} \pi(T)\pi(T')n^2\,.\end{equation}

Recall that the probability distribution $(\pi(T); T \subset [k])$ is derived from $\rho$ from equation (\ref{eqn:PIEpi}). To prove the equality above we begin by considering the ordered partition $\Pi$ associated to any $k$-tuple of subsets $(I^1, \ldots, I^k)$. The partition $\Pi$ has $2^k$ ordered cells $(\Pi(T); T \subset [k])$ defined by
$$\Pi(T) = \left(\bigcap_{t \in T} I^t \right) \cap \left(\bigcap_{t \notin T} [n]\setminus I^t \right ).$$
It follows from the inclusion-exclusion principle that if $(I^1, \ldots, I^k)$ has the density profile $\rho$ then $|\Pi(T)| = \pi(T)n$. The point here is that since $\pi$ can be derived from $\rho$, it in fact does not depend any individual $\Pi$.

For any fixed $k$-tuple $(I^1,\ldots,I^k)$, the sum $\sum_{\{u,v\}} \mathbf{1}_{\{T_u \cap T_v \neq \emptyset \}}$ can be represented by accounting for the contribution of each pair of subsets $\{T,T'\}$ to it.
$$\sum_{\{u,v\}} \mathbf{1}_{\{T_u \cap T_v \neq \emptyset \}} = \sum_{\{T,T'\}: T \cap T' \neq \emptyset} \# \left\{ \{u,v\}: T_u = T, T_v = T' \right\}.$$

Observe that by design $\Pi(T_u)$ is the cell of $\Pi$ that contains $u$, that is, $T_u = T$ if and only if $u \in \Pi(T)$. Therefore,
\begin{eqnarray*}
\# \left\{ \{u,v\}: T_u = T, T_v = T' \right\} &=& \# \left\{ \{u,v\}: u \in \Pi(T), v =\in \Pi(T') \right\} \\
&=& |\Pi(T)|\,|\Pi(T')| - \binom{|\Pi(T)\cap \Pi(T')|}{2} - |\Pi(T)\cap \Pi(T')|.\end{eqnarray*}

Since $|\Pi(T)| = \pi(T)n$ and $\Pi(T)\cap \Pi(T') = \emptyset$ for $T \neq T'$, the equality in (\ref{eqn:intersectioncount}) follows.
(The factor of 1/2 appears in (\ref{eqn:intersectioncount}) because we sum over all ordered pairs $(T,T')$.) Thus,
$$\E{Z(\rho)} \leq \# \{ k\text{\small{--tuples with density profile}}\; \rho \} (1-\frac{\lm}{n})^{\sum_{T \neq \emptyset} \binom{\pi(T)n}{2} + \frac{1}{2} \sum_{\substack{(T,T'): \\ T \neq T', T\cap T' \neq \emptyset}} \pi(T)\pi(T')n^2}.$$

The bijection between $k$-tuples and ordered partitions implies that the number of $k$-tuples with density profile $\rho$ is equal to the number of ordered partitions $(\Pi(T); T \subset [k])$ of $[n]$ such that $|\Pi(T) = \pi(T)n$. The latter number is the multinomial coefficient $\binom{n}{\pi(T)n\,; T \subset [k]}$. The statement of the lemma now follows.
\end{proof}

To analyze the asymptotic behaviour of $\E{Z(\rho)}$ write $\pi(T) = \beta(T) \frac{\log \lm}{\lm}$ for $T \neq \emptyset$. Then,
\begin{eqnarray*}
\sum_{T \neq \emptyset} \binom{\pi(T)n}{2} + \frac{1}{2} \sum_{\substack{(T,T') \\ T \neq T', T\cap T' \neq \emptyset}} \pi(T)\pi(T')n^2 &=& \frac{n^2}{2} \sum_{\substack{(T,T') \\ T\cap T' \neq \emptyset}} \pi(T) \pi(T') - \frac{n}{2} \sum_{T \neq \emptyset} \pi(T)
\end{eqnarray*}
and in terms of $\beta$ this equals 
\[\frac{n^2\log^2 \lm}{2\lm^2} \sum_{\substack{(T,T') \\ T\cap T' \neq \emptyset}} \beta(T) \beta(T') - \frac{n}{2} (1-\pi(\emptyset))\,.\]
Also, considering only the nonzero $\pi(T)$ and using Sterling's approximation we have
$$ \binom{n}{\pi(T)n\,;T \neq \emptyset} = O \left (n^{\frac{1-2^k}{2}} (\prod_{T} \pi(T))^{-1/2} \right ) \exp{\{nH(\pi)\}}$$
where $H$ is the previously introduced entropy function.

Using the fact that $1-\frac{\lm}{n} \leq e^{-\frac{\lm}{n}}$, and $1- \pi(\emptyset) \leq 1$, we conclude that
$$\E{Z(\rho)} = O \left ( n^{\frac{1-2^k}{2}} (\prod_{T} \pi(T))^{-1/2} e^{\lm/2} \right ) \exp{\left \{n \left[ H(\pi) - \frac{\log^2 \lm}{2\lm} \sum_{\substack{(T,T') \\ T\cap T' \neq \emptyset}} \beta(T) \beta(T') \right]\right \}}\,.$$

Recall in Lemma \ref{lem:asymptotics} we showed that $H(\pi) = \frac{\log^2 \lm}{\lm} \sum_{T \neq \emptyset} \beta(T) + O_k(\frac{\log \lm}{\lm})$, where the big O constant may depend on $k$. Consequently,
\begin{eqnarray} \label{eqn:pgwentropy}
H(\pi) - \frac{\log^2 \lm}{2\lm} \sum_{\substack{(T,T') \\ T\cap T' \neq \emptyset}} \beta(T) \beta(T') &=& \frac{\log^2 \lm}{\lm}\left(\sum_{T \neq \emptyset} \beta(T) - \frac{1}{2} \sum_{\substack{(T,T') \\ T\cap T' \neq \emptyset}} \beta(T) \beta(T')\right) \\
&+& O_k(\frac{\log \lm}{\lm}) \nonumber.
\end{eqnarray}

We also showed in Lemma \ref{lem:asymptotics} that if $\rho(S) = \alpha_{|S|} \frac{\log \lm}{\lm}$ for $S \neq \emptyset$ ($\rho(\emptyset) = 1$), then
\begin{equation} \label{eqn:alphabeta}
\sum_{T \neq \emptyset} \beta(T) - \frac{1}{2} \sum_{\substack{(T,T') \\ T\cap T' \neq \emptyset}} \beta(T) \beta(T') = \frac{1}{2} \sum_{i=1}^k (-1)^{i-1}\binom{k}{i} (2-\alpha_i)\alpha_i\,.
\end{equation}

Recall the event $A(\eps,p)$: the graph $G^i$ contains an independent set $S^i$ such that
the density profile of $(S^1, \ldots, S^k)$ satisfies
$$\rho(T) \in [(1-\eps) \alpha_{|T|,n,\lm,p}\frac{\log \lm}{\lm}, (1+\eps) \alpha_{|T|, n, \lm,p}\frac{\log \lm}{\lm}] \;\; \text{for all}\; T \subset [k].$$
We employ a first moment bound along with Lemma \ref{lem:pgwexpectation} to bound $\pr{A(\eps,p)}$; simplifying via (\ref{eqn:pgwentropy})
and (\ref{eqn:alphabeta}) we get
$$\pr{A(\eps,p)} \leq  \exp{\left \{ n\,[\frac{\log^2 \lm}{\lm} \Big (\sum_{i} (-1)^{i-1}\binom{k}{i} (\alpha_{i,n,\lm,p} - \frac{1}{2}\alpha^2_{i,n,\lm,p})
+ \mathrm{err}(\eps) \Big) + O_k(\frac{\log \lm}{\lm})] \right \}}$$
where $\mathrm{err}(\eps) \to 0$ as $\eps \to 0$ uniformly in $n, \lm$ and $p$.

From this point onward the proof of Theorem \ref{thm:pgwbinomeq} is completed in the same
manner as for regular graphs, which is the argument from Section \ref{sec:regulargraphconclusion}.

\subsection{A lower bound from regular trees} \label{sec:pgwlwbound}

We will show that factor of i.i.d.~independent sets in regular trees can be used to construct such independent sets in PGW trees as well.
Let $I$ be a factor of i.i.d.~independent set in the regular tree $\T_d$ such that the factor is a function of the labels in a finite size
neighbourhood of the root. Let $E(\lambda,d)$ denote the event that the root of $\pgw$ and all of its neighbours have degree at most $d$.

\begin{thm} \label{thm:coupling}
Given $I$ as above there exists a factor of i.i.d.~independent set $J$ of $\pgw$ whose density satisfies the bound
$$ \text{density}(I) \,\pr{E(\lambda,d)} \leq \text{density}(J) \leq \text{density}(I)\,.$$
\end{thm}

\begin{proof}
We construct $J$ in three stages.

\begin{description}
\item[The edge removal stage]  We remove edges from $\pgw$ via a factor of i.i.d.~process such that all vertices will have degree at most $d$ after the removal procedure. Begin with a random labelling $X$ of $\pgw$. For each vertex $v$ consider all the neighbours $u$ of $v$ such that the variables $X_u$ are the $\mathrm{degree}(v) - d$ highest in value (provided, of course, that $\mathrm{degree}(v) > d$). Mark the $\mathrm{degree}(v) - d$ edges connecting $v$ to these neighbours.

Following the marking procedure remove all the edges that have been marked. After the removal of edges, all vertices have degree at most $d$. The remaining graph is a disjoint collection of trees with a countable number of components. Denote it $G$.

\item[The filling out stage] If a vertex $v$ in $G$ has degree $\mathrm{degree}_G(v) < d$, then attach to it $d - \mathrm{degree}_G(v)$ copies of a $(d-1)$-ary tree via $d-1$ separate edges connecting $v$ to these trees. Following this procedure the graph $G$ becomes a disjoint collection of $d$-regular trees. Randomly label $G$ by a new set of labels $Y$ that are independent of $X$.

\item[The inclusion stage] Since $G$ is a disjoint collection of $d$-regular trees, we can use the factor associated to $I$ with input $Y$ to construct an independent set $I'$ of $G$ with the same density as $I$. Although $I'$ is an independent set in $G$ it may not be an independent set in the original tree $\pgw$ due to the removal of edges. To construct the independent set $J$, we include in $J$ all vertices $v \in I'$ such that no edges incident to $v$ were removed during the edge removal stage.
\end{description}

By design the random subset $J$ is a factor of i.i.d.~process on $\pgw$. $J$ is also an independent set because if $(u,v)$ is an edge of $\pgw$ with both $u, v \in I'$, then the edge connecting $u$ and $v$ must have been removed during the edge removal stage (due to $I'$ being an independent set in $G$). Thus neither $u$ nor $v$ belong to $J$.

To bound the density of $J$ we note that $J \subset I'$. Also, for any $v \in I'$, if $v$ and all of its neighbours in $\pgw$ has degree at most $d$ then none of the edges incident to $v$ are removed during the edge removal stage. Consequently, $v$ will be included in $J$. These two observations readily imply that
\[density(I) \,\pr{E(\lambda,d)} \leq density(J) \leq density(I)\,.\qedhere\]
\end{proof}
\begin{lem} \label{lem:poissontail}
If $\lambda = d - d^u$ for any $1/2 < u < 1$ then the probability $\pr{E(\lambda, d)} \to 1$ as $d \to \infty$.
\end{lem}
\begin{proof}
This is a calculation involving Poisson tail probabilities. Recall that the moment generating function of a $\poi(\mu)$ random variable
is $e^{\mu(e^t-1)}$. Let $X$ denote the degree of the root in a PGW tree of expected degree $\lambda$. Let $Z_1, \ldots, Z_X$
denote the number of offsprings of the neighbours of the root. Recall that $X$ has distribution $\poi(\lambda)$, and that conditioned on $X$ the random variables $Z_1, \ldots, Z_X$ are i.i.d. with distribution $\poi(\lambda)$.

Let $p(\lambda, d) = \pr{\poi(\lambda) > d}$. Then
\begin{eqnarray*}
\pr{E(\lambda,d)} &=& \E{\mathbf{1}_{X \leq d} \prod_{i=1}^X \mathbf{1}_{Z_i \leq d-1}}\\
&=& \E{\mathbf{1}_{X \leq d} \E{\prod_{i=1}^X \mathbf{1}_{Z_i \leq d-1} \mid X}}\\
&=& \E{\mathbf{1}_{X \leq d} (1-p(\lambda,d-1))^X}\\
&=& \E{(1-p(\lambda,d-1))^X} - \E{\mathbf{1}_{X > d}(1-p(\lambda,d-1))^X}\\
&\geq& e^{-\lambda \, p(\lambda,d-1)} - p(\lambda,d-1)
\end{eqnarray*}

We can bound the tail probability $p(\lambda,d-1)$ by using the exponential moment method. For simplicity we replace $d-1$ by $d$, which makes no difference to the analysis for large $d$. A simple and well-known computation gives
$$p(\lambda,d) \leq e^{d - \lambda}(\frac{\lambda}{d})^d \quad \text{if}\; \lambda < d \,.$$
Setting $\lambda = d - d^u$ for $1/2 < u < 1$, we see from the bound above that $p(d-d^u,d) \leq e^{d^u}(1-d^{u-1})^d = e^{d^u + d\log(1-d^{u-1})}$.
Since $\log(1-x) = - \sum_{k \geq 1} \frac{x^k}{k} \leq -x -x^2/2$ for $0 \leq x < 1$, by setting $x = d^{u-1} < 1$ we conclude that
\[ p(d-d^u,d) \leq e^{d^u - d(d^{u-1} + \frac{d^{2u-2}}{2})} = e^{-\frac{d^{2u-1}}{2}}\,.\]
Due to $u > 1/2$ the latter quantity tends to 0 exponentially fast as $d \to \infty$. As a result, both $(d-d^u)p(d-d^u,d)$ and $p(d-d^u,d)$ tend to 0 with $d$. This implies the lemma.
\end{proof}

Given $\lambda$, set $d = \lceil \lambda + \lambda^{3/4} \rceil$. From the definition of $\alpha(\lm), \alpha_d$, and
the conclusion of Theorem \ref{thm:coupling} we have that
$$\alpha(\lm) \frac{\log \lm}{\lm} \geq \alpha_d \cdot \frac{\log d}{d} \cdot \pr{E(\lm,d)}\,.$$

Recall the construction of Lauer and Wormald mentioned in Section \ref{sec:regulartrees} which shows that
$\liminf_{d \to \infty} \alpha_d \geq 1$. By our choice to $d$ as a function of $\lm$
we have $(\frac{\log d}{d})/( \frac{\log \lm}{\lm}) \to 1$ as $\lm \to \infty$. By Lemma (\ref{lem:poissontail})
we have that the probability $\pr{E(\lm, d)} \to 1$ as $\lambda \to \infty$.
As a result we conclude from the inequality above that
$$\liminf_{\lm \to \infty} \alpha(\lm) \geq 1\,.$$

This lower bound completes the proof of Theorem \ref{thm:pgwthm}.

\section{Concluding remarks} \label{sec:conclusion}

Our results are concerned with density of factor of i.i.d.~independent sets in sparse graphs where the sparsity parameter (degree) tends to infinity. However, it is still a very interesting problem to compute the maximum density for various classes of factor of i.i.d.~processes on $d$-regular graphs for fixed values of $d$. There have been some recent progress is this regard for independent sets in $3$-regular graphs. In \cite{CGHV} the authors use Gaussian processes to construct factor of i.i.d.~independent sets in $\T_3$ of density at least $0.436$, and in \cite{HW} the authors improve the bound to at least $0.437$ via another local algorithm . It is known due to McKay \cite{McKay} that the maximum density of independent sets in random 3-regular graphs is at most $0.456$, which provides an upper bound for factor of i.i.d.~independent sets in $\T_3$.

Another question is whether there is a gap between the density of the maximum cut in random $d$-regular graphs and the maximal density of
cuts derived from local algorithms. The density of the maximum cut of $\G_{n,d}$, denoted $\gamma(\G_{n,d})$, is the largest values
of $|\partial_{E}S|/nd$ where $\partial_{E}S$ is the set of all edges $(u,v) \in E(\G_{n,d})$ with $u \in S$ and $v \notin S$.
The expectation $\E{\gamma(\G_{n,d})} \to \gamma(d)$ as $n \to \infty$ for every $d$ \cite{BGT}. A local cut of $\T_d$ is a factor of i.i.d.~process
$\sigma \in \{-1,+1\}^{\T_d}$; its density is $\pr{\sigma(\circ) \neq \sigma(\circ')}$, where $(\circ,\circ')$ is a fixed edge of $\T_d$.
Is it true that $\gamma(d)$ equals the supremum over the density of local cuts of $\T_d$?

\end{document}